\documentclass{emsart}

\usepackage{stmaryrd,makeidx}

\contact[fwitt@math.fu-berlin.de]{Frederik Witt, FB Mathematik und Informatik, Freie Universit\"at Berlin, Arnimallee~3, D--14195 Berlin, FRG}

\newtheorem{theorem}{Theorem}[section]
\newtheorem{corollary}[theorem]{Corollary}
\newtheorem{lemma}[theorem]{Lemma}
\newtheorem{proposition}[theorem]{Proposition}

\newtheorem*{coro}{Corollary} 

\theoremstyle{definition}
\newtheorem{definition}[theorem]{Definition}
\newtheorem{remark}[theorem]{Remark}

\newenvironment{example}{\begin{trivlist}\item[]{\bf Example:}\setlength{\parindent}{0pt}}
{\end{trivlist}}

\def\note#1{\marginpar{\raggedright\if@twoside\ifodd\c@page\raggedleft\fi\fi\sf\scriptsize RMK: #1}}
\def\cliff{\mbox{\sl Cliff}}
\def\diff{\mbox{\sl Diff}}

\newcommand{\End}{\operatorname{End}}

\newcommand{\ric}{\operatorname{Ric}}
\newcommand{\eunorm}[1]{\parallel\mbox{\hspace{-3pt}} #1\mbox{\hspace{-3pt}}\parallel}

\newcommand{\htimes}{\hat{\otimes}}
\newcommand{\R}{\mathbb{R}}
\newcommand{\C}{\mathbb{C}}
\renewcommand{\H}{\mathbb{H}}
\newcommand{\Z}{\mathbb{Z}}

\newcommand{\mf}{\mathfrak}
\newcommand{\mb}{\mathbf}
\newcommand{\mc}{\mathcal}

\newcommand{\be}[1]{\begin{equation}\label{#1}}
\newcommand{\ee}{\end{equation}}
\newcommand{\gtilde}{\widetilde{\mc{G}}}

\title[Metric bundles of split signature and type II supergravity]{Metric bundles of split signature and type II supergravity}

\author[Frederik Witt]{Frederik Witt
\thanks{The author is supported by the DFG as a member of the SFB 647. He wishes to acknowledge the hospitality and the inspiring environment of the Erwin--Schr\"odinger--Institut, Wien. Moreover, he is also grateful to his collaborators Claus Jeschek and Florian Gmeiner (MPI f\"ur Physik, M\"unchen).}}

\begin{document}

\begin{abstract}
We study the geometry of type II supergravity compactifications in terms of an oriented vector bundle $E$, endowed with a bundle metric of split signature and further datum. The geometric structure is associated with a so--called generalised $G$--structure and characterised by an $E$--spinor $\rho$, which we can regard as a differential form of mixed degree. This enables us to reformulate the field equations of type II supergravity as an integrability condition of type $d_H\rho=0$, where $d_H=d+H\wedge$ is the twisted differential on forms. Finally, we investigate some geometric properties of integrable structures and formulate various no--go theorems.
\end{abstract}



\maketitle

\contentsline {section}{\numberline {1}Introduction}{2}
\contentsline {section}{\numberline {2}G--structures in supergravity}{4}
\contentsline {subsection}{\numberline {2.1}G--structures}{4}
\contentsline {paragraph}{Vector bundles.}{4}
\contentsline {paragraph}{Connections.}{7}
\contentsline {subsection}{\numberline {2.2}Clifford Algebras and spin structures}{8}
\contentsline {paragraph}{Clifford algebras.}{9}
\contentsline {paragraph}{Spin structures.}{12}
\contentsline {subsection}{\numberline {2.3}Compactification in supergravity}{14}
\contentsline {paragraph}{Heterotic supergravity.}{14}
\contentsline {paragraph}{Type II supergravity.}{16}
\contentsline {section}{\numberline {3}The linear algebra of $Spin(n,n)$}{17}
\contentsline {subsection}{\numberline {3.1}The group $Spin(n,n)$}{17}
\contentsline {subsection}{\numberline {3.2}Special orbits}{18}
\contentsline {paragraph}{Generalised $SU(m)$--structures\nobreakspace {}\cite {gmwi06},\nobreakspace {}\cite {jewi05a}.}{22}
\contentsline {paragraph}{Generalised $G_2$--structures\nobreakspace {}\cite {wi06}.}{22}
\contentsline {paragraph}{Generalised $Spin(7)$--structures\nobreakspace {}\cite {wi06}.}{23}
\contentsline {section}{\numberline {4}The generalised tangent bundle}{23}
\contentsline {subsection}{\numberline {4.1}Twisting with an $H$--flux}{23}
\contentsline {subsection}{\numberline {4.2}Spinors}{26}
\contentsline {section}{\numberline {5}The field equations}{30}
\contentsline {subsection}{\numberline {5.1}Integrable generalised $G$--structures}{30}
\contentsline {subsection}{\numberline {5.2}Geometric properties}{34}
%
%
%
%
%
\section{Introduction}

In this article, we will investigate the geometry of compactified type II supergravity theories by taking a $G$--structure point of view. 

Roughly speaking, a $G$--structure, for a given Lie group $G$, encodes the algebraic datum of a vector bundle -- the group is determined by looking at the symmetry group of a fibre inside its general linear group. For instance, a Riemannian metric and an orientation on the tangent bundle (i.e. an oriented Riemannian manifold $(M^n,g)$) is equivalent to an $SO(n)$--structure: Each tangent space defines an oriented Euclidean vector space $(T_pM,g_p)$ whose symmetry group inside $GL(T_pM)\cong GL(n)$ (i.e. the stabiliser under the $GL(n)$--action) is $SO(T_pM,g_p)\cong SO(n)$. The basic $G$--structure on $M^n$ underlying the subsequent development is an $SO(n,n)$--structure induced by an oriented vector bundle which carries a bundle metric of (``split") signature $(n,n)$. 

In a way, a $G$--structure is a device of putting a given  ``flat" model, i.e. a vector space for which an orientation, a metric or a complex structure etc. have been chosen, onto a manifold, where this datum varies smoothly from point to point. Regarding a vector space as a manifold, this datum is constant in the sense that the local coefficients with respect to the natural coordinate system are constant. In general however, we get non--trivial coefficients and the defect to be constant is measured by {\em integrability conditions}. For instance, on a Riemannian manifold $(M,g)$ there exists a coordinate system $(x_1,\ldots,x_n)$ such that $g=\sum dx_k^{2}$, if and only if the metric is flat, that is, the Riemannian curvature tensor vanishes. The idea of a $G$--structure approach to physical problems is then to interpret the {\em field content}, defined in physical space, for instance a gravitational field in general relativity, as the algebraic datum associated with a $G$--structure over some manifold, in this case a metric of signature $(1,3)$. The {\em field equations}, like vanishing Ricci tensor, are regarded as an integrability condition on the $G$--structure defined by the field content. This often simplifies the problem of finding solutions to the field equations dramatically, as representation theoretic arguments can be invoked to boil down the problem to mere linear algebra. We will illustrate this point and also recapitulate some general elements of vector bundle and representation theory as we go along. 

Before this, let us briefly introduce the mathematical content of the physical theory we propose to discuss. Neglecting the so--called {\em Ramond--Ramond fields}\index{Ramond--Ramond field} for the moment, the algebraic datum consists of a spinnable metric $g$, a closed integral 3--form $H$ (the {\em H--flux}), a scalar function $\phi$ (the {\em dilaton field}\index{dilaton field}) and 
two unit spinor fields $\Psi_{L,R}$ (the {\em supersymmetry parameters\index{supersymmetry parameters}}) -- we review spinors in some detail in Section~\ref{spinstruc}. In order to preserve two global supersymmetries (whence ``type II"\index{type II}), this datum is supposed to satisfy the following two supersymmetry equations\index{supersymmetry equations}. Firstly, there is the {\em gravitino equation}\index{gravitino equation}
\begin{equation}\label{susyeq1}
\nabla^{LC}_X\Psi_{L,R}\pm\frac{1}{4}(X\llcorner H)\cdot\Psi_{L,R}=0,
\end{equation}
where $\nabla^{LC}$ denotes the Levi--Civita connection of the metric and $X\llcorner H$ the $2$--form obtained from $H$ by contraction (inner product) with the vector field $X$. Secondly, the {\em dilatino equation} 
\begin{equation}\label{susyeq2}
(d\phi\pm\frac{1}{2}H)\cdot\Psi_{L,R}=0
\end{equation}
is supposed to hold. This situation is akin to so--called {\em heterotic}\index{heterotic} supergravity which has one global supersymmetry\index{supersymmetry} and is defined (in absence of additional gauge fields) by the gravitino\index{gravitino equation} and the dilatino equation\index{dilatino equation} on $\Psi_L$ alone. Here, $G$--structure techniques, where $G\subset Spin(n)$ is the stabiliser of a unit spinor, could be successfully applied to the construction of solutions (cf. for instance~\cite{ccdlmz03}, where the authors exploit the $SU(3)$\index{$SU(3)$}--structure associated with a unit spinor in dimension $6$). Similar attempts have been made for type II\index{type II} supergravity (see the review~\cite{gmw04}), but these turned out to be too restrictive. In a nutshell, the problem stems from the fact that if the pair $(\Psi_L,\Psi_R)$ is invariant under some $G\subset Spin(n)$, so is their common angle. From a physics point of view, however, the spinors $\Psi_L$ and $\Psi_R$ are independent, so the ``classical" $G$--structure ansatz coming from heterotic\index{heterotic} theory can only capture fairly particular cases. Rather, we should consider two independent $G$--structures associated with $\Psi_L$ and $\Psi_R$ simultaneously (for instance, two $SU(3)$\index{$SU(3)$}--structures if we deal with dimension $6$).

In~\cite{jewi05a} and~\cite{wi05}, it was noticed that the geometry of the supersymmetry equations\index{supersymmetry equations}~(\ref{susyeq1}) and~(\ref{susyeq2}) is captured by a so--called {\em generalised geometry}, a notion going back to the seminal article~\cite{hi03} of Hitchin. In mathematics, these structures naturally arose in connection with moduli spaces of geometric structures occurring in physics, where some points in the moduli space corresponded to {\em $B$--field transformations}\index{$B$--field transformation} of some well--known ``classical'' structure (cf. also~\cite{hu05}). Let us start by explaining what we understand by ``classical'' as opposed to ``generalised'' structures. Classically, any geometric structure on the tangent bundle can be transformed by diffeomorphisms. For instance, if $f:M\to M$ is a diffeomorphism and if $g$ defines a positive definite bundle metric on $TM$, so does the pulled--back tensor $f^*g$, that is both $(M,g)$ and $(M,f^*g)$ are Riemannian manifolds. In string theory, physicists are also used to transform geometric structures with $B$--fields, that is, $2$--forms on $M$. We thus want to pass from the natural transformation group $\diff(M)$ acting on $TM$ to the larger transformation group $\diff(M)\ltimes\Omega^2(M)$. How can this action of $2$--forms be implemented? Considering $B$ as a map from $TM\to T^*M$, we define $\exp(B):TM\oplus T^*M\to TM\oplus T^*M$ on sections $X\oplus\xi$ of $TM\oplus T^*M$ by
$$
\exp(B)(X\oplus\xi)=X\oplus\xi+\frac{1}{2}X\llcorner B.
$$
We will show in the sequel that this transformation can be indeed regarded as an exponential map as suggested by the notation. On the other hand, diffeomorphisms act on $TM\oplus T^*M$ in the usual way. The $B$--field transform of a ``classical'' geometry is then obtained by first extending the classical geometry from $TM$ to $TM\oplus T^*M$ in some suitable sense, and then by applying a $B$--field transformation\index{$B$--field transformation}. For instance, a Riemannian metric $g$ on $TM$ induces a map
$$
\mc{G}_0=\left(\begin{array}{cc} 0 & g^{-1}\\g & 0\end{array}\right)
$$
(where the matrix is taken with respect to the splitting $TM\oplus T^*M$ and we consider $g$ as a map from $TM\to T^*M$). Applying a $B$--field transformation\index{$B$--field transformation} gives the $B$--field transformed Riemannian metric
$$
\mc{G}_B=e^{2B}\circ\mc{G}\circ e^{-2B}=\left(\begin{array}{cc} -g^{-1}B & g^{-1}\\ g-Bg^{-1}B & Bg^{-1}\end{array}\right),
$$
an expression familiar to physicists (cf.~\cite{ka04}). We refer to $\mc{G}_B$ as a {\em generalised Riemannian metric}. The ``embedded'' case of $\mc{G}_0$ where $B=0$ is referred to as a {\em straight} generalised Riemannian metric.

The geometric structure of type II\index{type II} supergravity compactifications in dimension $6$ involves $B$--field transformed $SU(3)$\index{$SU(3)$}--structures, whose basic setup we describe next. The ``doubled'' bundle $TM\oplus T^*M$ carries a natural orientation and an inner product of split signature, namely contraction. Therefore, $TM\oplus T^*M$ is associated with an $SO(n,n)$--structure which is always spinnable, so we can consider $TM\oplus T^*M$--spinor fields. These can be regarded, modulo a scalar function, as differential forms of even or odd parity on $M$. A well--known procedure in spin geometry associates with the bitensor $\Psi_L\otimes\Psi_R$ a differential form $[\Psi_L\otimes\Psi_R]$. In the setup of type II\index{type II} supergravity compactified to dimension $6$, where the two unit spinors $\Psi_L$ and $\Psi_R$ induce an $SU(3)$\index{$SU(3)$}--structure each, one can show this form, being thought of as an $TM\oplus T^*M$--spinor, to be invariant under $SU(3)_L\times SU(3)_R\subset Spin(6,6)$\index{$SU(3)$} at every point. The datum $(M^6,g,\Psi_L,\Psi_R)$ induces therefore an $SU(3)_L\times SU(3)_R$\index{$SU(3)$}-- or {\em generalised $SU(3)$}\index{$SU(3)$}--structure on $TM\oplus T^*M$. This class of structures comprises the $B$--field transformation\index{$B$--field transformation} of ``classical'' $SU(3)$\index{$SU(3)$}--structures defined by $(M^6,g,\Psi)$. The straight generalised counterpart is defined by $[\Psi\otimes\Psi]$, that is, $\Psi_L=\Psi=\Psi_R$. The $B$--field transformation\index{$B$--field transformation} of a straight $SU(3)$\index{$SU(3)$}--structure is $\exp(B)\bullet[\Psi_L\otimes\Psi_R]$, where $\exp(B)\bullet$ stands for the action $(1+B\wedge+B\wedge B/2+\ldots)\wedge$ on differential forms.

The implementation of the $H$--flux $H$ requires one further idea, namely the concept of a {\em generalised tangent bundle}\index{generalised tangent bundle} $\mb{E}$~\cite{hi03},~\cite{hi05}. Since $H$ is closed, we can locally write $H_{|U_a}=dB^{(a)}$ for $B^{(a)}\in\Omega^2(U_a)$. On intersections $U_a\cap U_b$, we can twist the transition functions of $TM\oplus T^*M$ with $\exp(B^{(a)}-B^{(b)})$ which results in new transition functions inducing $\mb{E}$. If $H$ is integral this means that we are twisting with a {\em gerbe}. The field content on $M^6$, that is $g$, $H$, $\Psi_L$ and $\Psi_R$, therefore induces a generalised $SU(3)$\index{$SU(3)$}--structure on $\mb{E}$. Moreover, the supersymmetry equations\index{supersymmetry equations}~(\ref{susyeq1}) and~(\ref{susyeq2}) are equivalent to
\begin{equation}\label{formint}
d_He^{-\phi}[\Psi_L\otimes\Psi_R]=0,\quad d_He^{-\phi}[\mc{A}(\Psi_L)\otimes\Psi_R]=0,
\end{equation}
where $d_H=d+H\wedge$ is the {\em twisted differential} and $\mc{A}$ the {\em charge conjugation operator}. This is a natural integrability condition on the generalised $SU(3)$\index{$SU(3)$}--structure which links the supersymmetry equations\index{supersymmetry equations} to Hitchin's generalised variational principle~\cite{hi03},~\cite{wi06},~\cite{jewi05b}. En passant, we observe that this formalism is perfectly general and works for all classical $G$--structures defined by a spinor.

The spinorial formulation of integrability can be used to compute the Ricci and the scalar curvature of the metric, namely
$$
\ric(X,Y)=-2\mc{H}^{\phi}(X,Y)+\frac{1}{2}g(X\llcorner H,
Y\llcorner H),\quad S=2\Delta
\phi+\frac{3}{2}\eunorm{H}^2,
$$
where $\mc{H}^{\phi}=(\nabla^{LC})^2$ denotes the Hessian of the dilaton field\index{dilaton field}, and $\Delta$ the Laplacian of the metric. Moreover, this description gives rise to two striking no--go theorems in the vein of similar statements from~\cite{gmpw04} and~\cite{wi06}. Firstly, integrability implies $H=0$ if $M$ is compact. This means that both spinors $\Psi_{L,R}$ are parallel with respect to the Levi--Civita connection, and so the holonomy reduces to the intersection of the two $G$--structures associated with $\Psi_{L,R}$ (this situation belongs indeed to the realm of classical $G$--structures, for the spinors are either parallel or may be assumed to be orthogonal by parallelity). Hence, there are no interesting generalised geometries satisfying~(\ref{formint}) for $M$ compact. Secondly, and independently of compactness, $H=0$ if and only if $d\phi=0$, so in order to obtain non--parallel solutions to~(\ref{susyeq1}) and~(\ref{susyeq2}), we need indeed a non--trivial dilaton field\index{dilaton field}. The construction of solutions with non--trivial $H$--flux is thus a fairly difficult task, and only local examples are known so far~\cite{wi06}.
%
%
%
%
%
\section{G--structures in supergravity}
\label{sugra}
In this section we first recall the basic features of $G$--structure theory and spin geometry as far as we need it here. We will then briefly outline the theory of heterotic\index{heterotic} and type II\index{type II} supergravity compactifications, thereby illustrating the use of $G$--structure techniques. 
%
%
\subsection{G--structures}
\label{gstructures}
Suitable references for this section are~\cite{kono63} and~\cite{mo01}.

\paragraph{Vector bundles.}
\vspace{-10pt}
Let $M^n$ be a differentiable manifold of dimension $n$ and $\pi:\mb{V}^m\to M^n$ be a real vector bundle of rank $m$. By definition, there exists an open cover $\{U_a\}$ of $M$ with trivialisations $s_a:U_a\times\R^m\stackrel{\cong}{\to}\mb{V}^m_{|U_a}$. If $x_1,\ldots,x_m$ denotes the standard basis of $\R^m$, $s_{a,k}(p)=s_a(p,x_k)$ defines a basis of the fibre $\mb{V}^m_p$. For $p\in U_a\cap U_b$, the bases $s_{a,k}(p)$ and $s_{b,k}(p)$ are related by an element of $GL(m)$, so we get a collection of {\em transition functions} $s_{ab}:U_a\cap U_b\to GL(m)$. These satisfy the {\em cocycle condition} 
\begin{equation}\label{cocycle}
s_{ab}\circ s_{bc}=s_{ac},
\end{equation} 
whenever $U_a\cap U_b\cap U_c\not=0$. 
Conversely, given a collection of functions $\{s_{ab}:U_a\cap U_b\to GL(m)\}$ such that~(\ref{cocycle}) holds, we can define a real, rank $m$ vector bundle $\mb{V}^m$ with the $GL(m)$--module $V^m$ as fibre
$$
\mb{V}^m=\coprod\limits_a U_a\times V^m/\sim_{s_{ab}}.
$$
Here, two elements $(a,p,v)$, $(b,q,w)$ are equivalent if and only if $p=q$ and $v=s_{ab}(w)$. The projection is defined by $\pi([a,p,v])=p$. Local trivialisations are provided by
$$
s_a:(p,v)\in U_a\times V^m\mapsto[p,a,v]\in\mb{V}^m,
$$
which over $U_a\cap U_b$ induce again the transition functions $s_{ab}$. In this picture, a {\em section} of $\mb{V}^m$, i.e. a smooth map $\sigma:M\to\mb{V}^m$ such that $\pi\big(\sigma(p)\big)=p$, is a collection of maps $\sigma=\{\sigma_a:U_a\to V^m\}$ with $\sigma_a=s_{ab}(\sigma_b)$, so that $\big(a,p,\sigma_a(p)\big)\sim\big(b,p,\sigma_b(p)\big)$ for $p\in U_a\cap U_b$. We denote the $C^{\infty}(M)$--module of sections of $\mb{V}^m$ by $\Gamma(\mb{V}^m)$.

Given the vector bundle $\mb{V}^m$, we can consider further bundles derived from $\mb{V}^m$, associated with $GL(m)$--modules derived from $V^m$. For instance, $V^{m*}$ gives the {\em dual bundle} $$
\mb{V}^{m*}=\coprod\limits_a U_a\times V^{m*}/\sim_{s_{ab}},
$$ 
$\Lambda^pV^{m*}$ the {\em bundle of exterior $p$--forms} 
$$
\Lambda^p\mb{V}^{m*}=\coprod\limits_a U_a\times \Lambda^pV^{m*}/\sim_{s_{ab}},
$$ 
and $\mf{gl}(m)$, the Lie algebra of $GL(m)$, the {\em adjoint bundle} 
$$
Ad(\mb{V}^m)=\coprod\limits_a U_a\times \mf{gl}(m)/\sim_{s_{ab}}.
$$ 
In fact, for any $G$--space $F$ we obtain the {\em associated fibre bundle} 
$$
\mathbb{F}=\coprod U_a\times F/\sim_{s_{ab}}
$$ 
with typical fibre $F$.

\begin{remark}
This local approach to vector and fibre bundles in terms of open covers is dissatisfactory insofar as this involves the choice of a specific collection of trivialisations which is far from being canonic in the same way as a different choices for an atlas of a differentiable manifold can be made. This can be circumvented by using principal $G$--fibre bundles, but from a practical point of view, the local description gives a suitable working definition for the later development and is closer to the physical intuition which is why we use it here.
\end{remark}

We say that a vector bundle $\pi:\mb{V}^m\to M^n$ carries a {\em $G$--structure}\index{$G$--structure}, if there exists an open cover $\{U_a\}$ of $M^n$ with trivialisations of $s_a:U_a\to \mb{V}_{|U_a}$, whose induced transition functions take values in $G$. We also speak of a {\em reduction}\index{reduction} of the structure group of $\mb{V}^m$ to $G$. {\em Linear} $G$--structures\index{$G$--structure, linear} are associated with the tangent bundle, so $G\subset GL(n)$. Two other structures will be important to us: {\em spin structures}\index{spin structure}, where $G\subset Spin(n)$ (cf. Section~\ref{spinstruc}), and {\em generalised structures}, where $G\subset SO(n,n)$ or $Spin(n,n)$ (cf. Section~\ref{linalg}).

The group $G$ acts on any $GL(m)$--representation space via restriction. To understand where $G$--structures come from, we first have to understand the underlying $G$--representation theory. We briefly recall some concepts we will make intensive use of later on (a good reference is~\cite{fuha91}). A {\em representation}\index{representation} of a group $\mb{G}$ consists of a vector space $V$ and a smooth group homomorphism $\rho_V:\mb{G}\to GL(V)$. The simplest example is the trivial representation where $\rho_V\equiv Id_V$. Therefore, $V$ becomes a $\mb{G}$--module under the action of $\mb{G}$, and in particular, we get a disjoint decomposition into orbits of the form $\mb{G}/G$ for some subgroup $G$. The $G$--structures we will consider in the sequel arise precisely in this way. In general, the determination of the orbit structure on $V$ is a difficult problem. One therefore rather looks for $\mb{G}$--invariant subspaces of $V$. If there is no invariant subspace other than $\{0\}$ and $V$ itself, the representation is said to be {\em irreducible}. For large classes of groups (for instance, $\mb{G}$ compact or semi--simple), the {\em complete reducibility property} holds: Any representation space can be decomposed into a direct sum of irreducible subspaces. A linear map $F:V\to W$ between $\mb{G}$--representation spaces is $\mb{G}$--{\em equivariant}\index{equivariant} if the action commutes with $F$, i.e. 
$$
F\big(\rho_V(g)(v)\big)=\rho_W(g)\big(F(v)\big).
$$ 
Two $\mb{G}$--representations $V$ and $W$ are {\em equivalent}\index{equivalent}, if there exists a $\mb{G}$--equivariant isomorphism. If both representation spaces $V,W$ are irreducible, Schur's Lemma\index{Schur's Lemma} asserts $F$ to be either an isomorphism or to be trivial. 

Given the way $G$ arises, our first task is to understand the quotient $GL(m)/G$ and in particular, its trivial coset $[Id]$. It represents a set of $G$--invariant objects $Q_1,\ldots,Q_r$.

\begin{example}
Take $G=O(m)\subset GL(m)$ which stabilises a Euclidean metric $Q_1=g$. The map 
$$
[A]\in GL(m)/O(m)\mapsto A^*g\in\odot^2_+\R^{m*},
$$ 
where $A\in GL(m)$ acts via $A^*$ in the standard way on symmetric tensors (in its own right an honest representation $GL(m)\to GL(\odot^2\R^{m*})!)$, sets up a bijection between the coset space $GL(m)/O(m)$ and the set of positive definite symmetric $2$--tensors over $\R^m$, under which $[Id]$ corresponds to $g$. There are preferred ``$G$--bases'' in which the invariants $Q_1,\ldots,Q_r$ take a special shape, in our example orthonormal bases $e_1,\ldots,e_m$ with dual basis $e^1,\ldots,e^m$, for which $g=\sum e^k\otimes e^k$. Passing to global issues, assume that the vector bundle $\mb{V}^m$ admits at least one collection of $O(m)$--valued transition functions $\{s_{ab}\}$, i.e. $\mb{V}^m$ carries an $O(m)$--structure. We define local sections of the derived bundle $\odot^2\mb{V}^{m*}$ by 
\begin{equation}\label{metricten}
g_a=\sum s_a(p,e^k)\otimes s_a(p,e^k).
\end{equation} 
Since the associated transition functions take values in $O(m)$, the local orthonormal basis $s_a(p,e_k)$ gets mapped to the local orthonormal basis $s_b(p,e_k)$ under $s_{ab}$. Consequently, the collection $\{g_a\}$ patches together to a global section of $\odot^2\mb{V}^{m*}$ and thus defines a bundle metric on $\mb{V}^m$. Conversely, assume to be given a bundle metric $g$. Then $g$ singles out preferred local bases, namely those for which $g_a=g_{|U_a}$ has standard form~(\ref{metricten}). The resulting transition functions take values in $O(m)$. In particular, a Riemannian manifold $(M^n,g)$ is nothing else than an $O(n)$--structure. Note that for an arbitrary local trivialisation of $\mb{V}^{m*}$, $g_{|U_a}$ will not acquire its standard form. Nevertheless, in order to define the metric globally, it is enough to exhibit one cover for which it does. 
\end{example}

In the same way, the invariants $Q_1,\ldots,Q_r$ acquire global meaning for an arbitrary $G$--structure. Besides the $G$--invariants, the decomposition of $V^m=\oplus_kV_k$ into irreducible $G$--modules also carries over to a global decomposition of $\mb{V}^m=\oplus_k\mb{V}_k$ where $\mb{V}_k$ is the vector bundle with fibre $V_k$.

\begin{example}
For $G=O(m)$, the vector representation is of course irreducible, but $\odot^2V^{m*}=\mb{1}\oplus\odot_0^2V^{m*}$ (the second summand consisting of the trace--free symmetric endomorphisms of $V^m$). Hence we obtain an analogous decomposition of $\odot^2\mb{V}^{m*}=\mb{1}\oplus\odot_0\mb{V}^{m*}$. 
\end{example}

Within a particular $G$--structure, further reductions are possible. For instance, we can reduce $O(m)$ to $SO(m)$. The coset space $O(m)/SO(m)$ is isomorphic to $\Z_2$ and can be identified with the set $\Lambda^m\R^{m*}/\R_{>0}$ (i.e. two volume forms being equivalent if they differ by a positive real scalar) by sending $[A]$ to the class $[A^*(e^1\wedge\ldots\wedge e^m)]=[\det(A^{-1})\cdot e^1\wedge\ldots\wedge e^m]$. In analogy with the previous example, this yields a globally defined volume form and therefore an orientation on each fibre of $\mb{V}^m$. Consequently, the real line bundle $\Lambda^m\mb{V}^{m*}$ is trivial, which forces its characteristic class, namely the {\em first Stiefel--Whitney class} $w_1(\mb{V}^m)$, to vanish. We thus meet a topological obstruction against reducing the structure group to $SO(n)$. In general, reductions from $\mb{G}$ to $G$ are parametrised by sections of the fibre bundle $\coprod U_a\times\mb{G}/G/\sim_{s_{ab}}$. Since the existence of such sections is a purely topological question, as highlighted in the previous example, one also speaks of {\em topological} reductions\index{reduction, topological} as opposed to {\em geometrical} reductions\index{reduction, geometrical}. We turn to these next.

\paragraph{Connections.}
To do differential geometry over $\mb{V}^m$ requires the choice of a covariant derivative\index{covariant derivative}.

\begin{definition}\label{connection}
A {\em linear connection}\index{connection} or {\em covariant derivative}\index{covariant derivative} on a vector bundle $\mb{V}^m$ is a linear map
$$
\nabla:\Gamma(\mb{V}^m)\to\Gamma(T^*M\otimes\mb{V}^m),
$$
such that 
\begin{equation}\label{prodrule}
\nabla(f\sigma)=df\otimes\sigma+f\cdot\nabla\sigma
\end{equation} 
holds for any real function $f\in C^{\infty}(M)$. 
\end{definition}

We usually speak simply of a connection\index{connection} for short. Again, following our general philosophy, we want to characterise a connection in local terms. So we fix a local basis $s_{a,k}=s_a(\cdot,e_k)$ of $\mb{V}^m$ over $U_a$, where $e_1,\ldots,e_m$ is the standard basis of the fibre $\R^m$. If $\sigma\in\Gamma(\mb{V}^m)$, then $\sigma_{|U_a}=\sum f^k_as_{a,k}$ for smooth functions $f^k_a\in C^{\infty}(U_a)$. By~(\ref{prodrule}),
$$
\nabla\sigma_{|U_a}=\sum df^k_a\otimes s_{a,k}+f^k_a\cdot\nabla s_{a,k}.
$$
Now $\nabla s_{a,k}$ is again a section of $\mb{V}^m_{|U_a}$, hence $\nabla s_{a,k}=\sum\omega_{a,k}^je_j$ for $1$--forms $\omega_{a,k}^j$ (i.e. $\nabla_Xs_{a,k}=\sum\omega_{a,k}^j(X)s_{a,j}$ for any vector field $X$ on $U_a$). Using matrices, the action of the connection is described by
$$
\nabla \left(\begin{array}{c}f^1_a\\ \vdots\\f^m_a\end{array}\right)=\left(\begin{array}{c}df^1_a\\ \vdots\\df^m_a\end{array}\right)+\left(\begin{array}{ccc}\omega^1_{a1} & \cdots & \omega^1_{am}\\ \vdots & \ddots & \vdots\\ \omega^m_{a1} & \ldots & \omega_{am}^m \end{array}\right)\left(\begin{array}{c}f^1_a\\\vdots\\f^m_a\end{array}\right)
$$
Locally, a connection is therefore just a differential operator of the form $d+\omega_a$ for some $m\times m$ matrix $\omega_a$ of $1$--forms defined over $U_a$. Since the connection\index{connection} is a global object, it is in particular defined on overlaps. Now the transformation rule $\tau_a=s_{ab}(\tau_b)$ for global section $\tau\in\Gamma(\mb{V}^m)$ thought of as a family $\{\tau_a:U_a\to\R^m\}$ implies for $\tau=\nabla\sigma$
\begin{equation}\nonumber
d\sigma_a+\omega_a\sigma_a=s_{ab}(d\sigma_b+\omega_b\sigma_b).
\end{equation}
A straightforward computation shows 
\begin{equation}\label{connectiontrans}
\omega_a=s_{ab}ds_{ab}^{-1}+s_{ab}\omega_bs_{ab}^{-1},
\end{equation}
implying the following local characterisation of a connection\index{connection}.

\begin{proposition}
A connection\index{connection} is given by a collection of smooth matrices of $1$--forms $\{\omega_a:U_a\to\R^{m\times m}\}$ such that~(\ref{connectiontrans}) holds.
\end{proposition}

To make contact with $G$--structures, we first note that the space of $m\times m$--matrices can be identified with the Lie algebra $\mathfrak{gl}(m)$ of $GL(m)$. If $\theta$ is the tautological $1$--form of $GL(m)$, $\theta_{ab}=s_{ab}^*\theta$ the pull--back to $U_a\cap U_b$ and $Ad$ denotes the adjoint action of $GL(m)$ on its Lie algebra, the gluing rule~(\ref{connectiontrans}) reads
\begin{equation}\nonumber
\omega_a=Ad(s_{ab})\omega_b+\theta_{ab}.
\end{equation}
The collection $\omega=\{\omega_a\}$ thus consists of local, $\mf{gl}(m)$--valued maps which is why these connections are referred to as ``linear'' (cf. Definition~\ref{connection}). As a result, a connection can be defined with respect to a given collection of $GL(m)$--valued transition functions without explicit reference to the vector bundle $\mb{V}^m$. In particular, a $GL(m)$--connection induces a $\nabla$--operator on any vector bundle associated with $GL(m)$--valued transition functions. For example, if we can covariantly derive vector fields (i.e. sections of the tangent bundle), we get a canonic covariant derivative\index{covariant derivative} for any tensor bundle.
 
Now given a fixed connection $\nabla$ and a topological reduction\index{reduction, topological} from a $GL(m)$-- to a $G$--structure, we refer to this reduction as {\em geometrical}\index{reduction, geometrical} if the $\omega_a$ take values in $\mf{g}$, the Lie algebra of $G$, rather than in $\mf{gl}(m)$. Since $\mf{g}$ acts trivially on the $G$--invariant objects $Q_1,\ldots,Q_r$ ($G$ acting as the identity), a connection reduces geometrically to $G$ if and only if $\nabla Q_k=0$ for $k=1,\ldots,r$ (where $\nabla$ is extended to the corresponding vector bundles). As geometrical reductions presuppose the choice of a covariant derivative\index{covariant derivative}, this notion is particularly interesting if we can make a canonic choice for $\nabla$.

\begin{example}
A connection is {\em metric} if and only if it reduces to a given $O(m)$--structure, that is, if and only if $\nabla g=0$. An example of this is the {\em Levi--Civita covariant derivative} $\nabla^{LC}$ of a Riemannian manifold $(M^n,g)$ (i.e. an $O(n)$--structure), which is implicitly defined by
\begin{eqnarray*}
2g(\nabla^{LC}_XY,Z) & = & \phantom{-}X.g(Y,Z)+Y.g(Z,X)-Z.g(X,Y)\\
&  & -g(X,[Y,Z])+g(Y,[Z,X])+g(Z,[X,Y]).
\end{eqnarray*}
It is the unique metric connection $\nabla$ whose {\em torsion tensor} 
$$
T^{\nabla}(X,Y)=\nabla_XY-\nabla_YX-[X,Y]
$$ 
vanishes. The Levi--Civita connection geometrically reduces geometrically to $G$ if and only if the holonomy group of $g$ is contained in $G$. The $G$--structures to which $\nabla^{LC}$ can reduce are therefore given by by Berger's famous list~\cite{be55}.
\end{example}
%
%
\subsection{Clifford Algebras and spin structures\index{spin structure}}
\label{spinstruc}

In relativistic particle physics, particles arise as elements of an irreducible representation for a symmetry group $G\subset Spin(1,q)$. They come in two flavours: They are either {\em bosonic}\index{bosonic} (particles that transmit forces) and are elements of a {\em vector representation}\index{vector representation} of $Spin(1,q)$, or they are fermionic\index{fermionic} (particles that make up matter) and live in a {\em spin representation}\index{spin representation} of $Spin(1,q)$, a notion we now review. For details on the aspects treated below, we recommend~\cite{ba81},~\cite{ha91},~\cite{lami89} and~\cite{wa89}.

\paragraph{Clifford algebras.}
Let $V$ be a real or complex vector space equipped with a metric $g$, i.e. a symmetric non--degenerate bilinear form. Out of the datum $(V,g)$ we can construct the {\em Clifford algebra}\index{Clifford algebra} $\cliff(V,g)$ as a deformation of the exterior algebra, namely

\begin{itemize}
\item as a vector space, $\cliff(V,g)=\Lambda^*V\stackrel{g}{\cong}\Lambda^*V^*$, and in particular, $V\subset\cliff(V,g)$.
\item with an algebra product subject to the relation $X\cdot Y+Y\cdot X=-2g(X,Y)\mb{1}$ for $X,\,Y\in V$.
\end{itemize}

In particular, this induces a grading into elements of degree $p$ and a coarser $\Z_2$--grading into even and odd elements $\cliff(V,g)^{ev,od}$. As any real (pseudo--) Euclidean vector space $(V,g)$ is isometric to some $\R^{p,q}$, the resulting real Clifford algebras\index{Clifford algebra} are -- up to isomorphism -- given by $\cliff(p,q)=\cliff(\R^{p,q})$. Here, the notation $\R^{p,q}$ refers to the vector space $\R^{p+q}$ endowed with its canonical inner product of signature $(p,q)$. We simply write $\cliff(p)$ for $\cliff(p,0)$. Moreover, $\cliff(V^{\C},g^{\C})\cong\cliff(V,g)\otimes\C$. Since the complexification of a metric $g$ depends only on the dimension of $\R^{p,q}$, but not on the signature, there is only one type of complex Clifford algebras\index{Clifford algebra}, depending on the dimension of $V^{\C}$. 

In universal terms, one can characterise $\cliff(V,g)$ (up to algebra isomorphism) as the unique algebra satisfying the following property: If $A$ is an associative algebra with unit over the same field as $V$, and $f:V\to A$ is a linear map such that $f(X)\cdot f(X)=-g(X,X)\mb{1}$, then $f$ extends to an algebra homomorphism $\cliff(V,g)\to A$ in a unique way. As a consequence, Clifford algebras\index{Clifford algebra} are essentially matrix algebras in disguise. 
To see this, choose an almost complex structure $J$ on $\R^{2m}$ which acts as an isometry for $g$, together with a real subspace $U$, defining an orthogonal splitting $\R^{2m}=U\oplus J(U)$. We endow the complexification $U^{\C}$ of $U$ with the hermitian inner product $q$ induced by $g_{|U}$ and extend it to $\Delta_{2m}=\Lambda^*U^{\C}$, the space of {\em Dirac spinors}\index{Dirac spinor}. For $u\in U^{\C}$, let $u\llcorner$ be the hermitian adjoint of $u\wedge$ with respect to $q$, which by convention we take to be conjugate--linear in the first argument. We identify $\R^{2m}$ with $U^{\C}$ via $j(u_1\oplus Ju_2)=u_1+iu_2$. Extension of the map 
$$
f:x\in\R^{2m}\mapsto f_x\in\End(\Delta_{2m}),\quad f_x(\psi)= j(x)\wedge\psi-j(x)\llcorner\psi
$$
to $\cliff(2m,0)\otimes\C$ squares to minus the identity for unit vectors. Hence, by the universal property, 
$$
\cliff(\C^{2m})\cong\End(\Delta_{2m}).
$$ 
The resulting action is usually referred to as {\em Clifford multiplication}\index{Clifford multiplication} and is denoted by a dot, so $a(\Psi)=a\cdot\Psi$ for $a\in\cliff(\C^{2m})$. With respect to the hermitian inner product $q$, we have for $a\in\cliff(\C^{2m})$, $\Psi,\,\Phi\in\Delta_{2m}$,
$$
q(a\cdot\Psi,\Phi)=q(\Psi,\widehat{a}\cdot\Phi),
$$
where $\wedge$ is a sign--changing operator defined on elements of degree $p$ to be 
$$
\widehat{a}=(-1)^{p(p+1)/2}a.
$$ 
For dimension $2m-1$, we consider as above the even--dimensional space $\R^{2m}=U'\oplus J(U')$ and put $\R^{2m-1}=U'\oplus J(U)$ where $U'=U\oplus\R e$ for some unit vector $e$. Since the even part of a Clifford algebra\index{Clifford algebra} is an algebra in its own right, the map $\gamma:\R^{2m-1}\to\cliff(2m,0)^{ev}$ defined by $\gamma(x)=x\cdot e$ extends to an algebra homomorphism, which actually is an isomorphism. Hence $\cliff(\C^{2m-1})$ acts on $\Delta_{2m}$ via restriction of the action of $\cliff(\C^{2m})$ to the image under $\gamma$. We are therefore left with determining a representation of $\cliff(\C^{2m})^{ev}$. Now for an orthonormal basis $e_1,\ldots,e_{2m}$, $\omega^{\C}=(-1)^{m(m+1)/2}i^me_1\wedge\ldots\wedge e_{2m}$ defines an involution on $\Delta_{2m}$ whose $\pm1$--eigenspaces $\Delta_{\pm}\subset\Delta$ contain the so--called {\em Weyl spinors}\index{Weyl spinor}. The subscript indicates the {\em chirality}\index{chirality} of the spinor. Chirality is preserved under the action of $\cliff(\C^{2m})^{ev}$ while it is reversed under odd elements, so
$$
\cliff(\C^{2m})^{ev}=\End(\Delta_+)\oplus\End(\Delta_-).
$$ 
As $\cliff(\C^{2m-1})$--representations, $\Delta_+$ and $\Delta_-$ are equivalent, i.e. there is an isomorphism commuting with the action of $\cliff(\C^{2m-1})$ such that $\Delta_+\cong\Delta_-\cong\Delta_{2m-1}$. Consequently, we find
$$
\cliff(\C^{2m+1})= \End_{\C}(\Delta_{2m-1})\oplus\End(\Delta_{2m-1}),
$$
from which we obtain a representation of $\cliff(\C^{2m-1})$ by projecting on the first factor (the choice being immaterial). The even part can then be identified with $\cliff(\C^{2m-1})^{ev}\cong\End(\Delta_{2m-1})$. As a vector space, $\Delta_{2m-1}$ is isomorphic to $\C^{2^{m-1}}$. Summarising, we obtain the $\mod2$--periodicity
$$
\cliff(\C^{2m})\cong\End_{\C}(\Delta_{2m}),\quad\cliff(\C^{2m+1})\cong\End_{\C}(\Delta_{2m+1})\oplus\End_{\C}(\Delta_{2m+1}).
$$
In the real case, a careful analysis reveals a mod $8$ periodicity (which, as the mod $2$ periodicity in the complex case, is ultimately an instantiation of Bott periodicity) depending on the signature $(p,q)$ of $g$, namely
\begin{equation}\label{clifftable}
\begin{array}{cl}
p-q\mod8: & \\
0,6 & \quad\cliff(p,q)\cong\End_{\R}(P_{p,q})\\
1,5 & \quad\cliff(p,q)\cong\End_{\C}(P_{p,q})\\
2,4 & \quad\cliff(p,q)\cong\End_{\H}(P_{p,q})\\
3   & \quad\cliff(p,q)\cong\End_{\H}(P_{p,q})\oplus\End_{\H}(P_{p,q})\\
7   & \quad\cliff(p,q)\cong\End_{\R}(P_{p,q})\oplus\End_{\R}(P_{p,q}),
\end{array}
\end{equation}
where $P_{p,q}$ is a vector space over the appropriate ground field which we view as a real vector space when necessary. For the even parts, we find
\begin{equation}\label{cliffevtable}
\begin{array}{cl}
p-q\mod8: & \\
0 & \quad\cliff(p,q)^{ev}\cong\End_{\R}(S_{p,q+})\oplus\End_{\R}(S_{p,q-})\\
1,7 & \quad\cliff(p,q)^{ev}\cong\End_{\R}(S_{p,q})\\
3,5 & \quad\cliff(p,q)^{ev}\cong\End_{\H}(S_{p,q})\\
2,6 & \quad\cliff(p,q)^{ev}\cong\End_{\C}(S_{p,q+})\cong\End_{\C}(S_{p,q-})\\
4   & \quad\cliff(p,q)^{ev}\cong\End_{\H}(S_{p,q+})\oplus\End_{\H}(S_{p,q-}),
\end{array}
\end{equation}
where $S_{p,q(\pm)}$ is a vector space over the appropriate ground field which again we view as a real vector space when necessary. As in the complex case, we refer to the elements of $S_{p,q\pm}$ as spinors of positive or negative chirality. The cases which concern us most are $p-q\mod8=0,6,7$; similar remarks apply to the remaining cases. For $p-q\mod8=0$, the $\cliff(p,q)$--representation space $P_{p,q}$ is real and can be decomposed into the $\pm1$--eigenspaces $S_{p,q\pm}$ of the Riemannian volume form $vol_g$ which defines an involution. As a result, we find $\Delta_{p+q\pm}=S_{p,q\pm}\otimes\C$. Similarly, the representation of $\cliff(p,q)$ with $p-q\mod 8=7$ is real, and the spin representation\index{spin representation} is just $S_{p,q}=P_{p,q}$, whence $\Delta_{p+q}=S_{p,q}\otimes\C$. More interesting is the case $p-q\mod8=6$. Here, the $\cliff(p,q)$--representation space is also real. However, as displayed in Table~(\ref{cliffevtable}), $\cliff(p,q)^{ev}$ is isomorphic with $\End_{\C}(S_{p,q+})$ considered as a real algebra, and therefore also to $\End_{\C}(\overline{S}_{p,q+})$, where $S_{p,q-}=\overline{S}_{p,q+}$ is the conjugated representation\footnote{If $\rho:\mb{G}\to GL_{\C}(V)$ is a representation over a complex vector space, then the {\em conjugate representation} is defined by $\overline{\rho}(g)=\overline{\rho(g)}$, i.e. we complex conjugate the entries of the matrix $\rho(g)$.}. Hence, $\Delta_{p+q}=P_{p,q}\otimes\C=S_{p,q+}\oplus\overline{S}_{p,q+}$, so that $S_{p,q\pm}=\Delta_{p+q\pm}$ are the $\pm i$--eigenspaces of $vol_g$.

The Clifford algebras\index{Clifford algebra} $\cliff(p,q)$ define various substructures of interest, for instance the {\em Pin} and the {\em Spin} group of a pseudo--Euclidean vector space $(\R^{p,q},g)$. Let
\begin{equation}\nonumber
Pin(p,q)=\{x_1\cdot\ldots\cdot x_l\,|\,x_j\in \R^{p,q},\,\eunorm{x_j}=\pm1\}\subset\cliff(p,q).
\end{equation}
The spin group $Spin(p,q)\subset\cliff(p,q)^{ev}$ is then the subgroup generated by elements of even degree. The covering map 
\begin{equation}\label{covering}
\pi_0:a\in Pin(p,q)\mapsto \pi_0(a)\in O(p,q),\qquad\pi_0(a)(x)=(-1)^{\deg(a)}a\cdot x\cdot a^{-1}
\end{equation}
gives rise to the exact sequences
$$
\begin{array}{l}
\{1\}\to\Z_2\to Pin(p,q)\stackrel{\pi_0}{\to} O(p,q)\to \{Id\}\\[5pt]
\{1\}\to\Z_2\to Spin(p,q)\stackrel{\pi_0}{\to} SO(p,q)\to \{Id\}.
\end{array}
$$
Note that $Pin(p,q)$ consists of four connected components unless $p=0$ or $q=0$ (in which case there remain only two). Similarly, $Spin(p,q)$ has two connected components for $p,q\geq1$. We denote by $Spin(p,q)_+$ the identity component of $Spin(p,q)$ which as a set is 
\begin{equation}\nonumber
Spin(p,q)_+=\{x_1\cdot\ldots\cdot x_{2l}\,|\,\eunorm{x_j}=1\mbox { for an even number of $j$}\}.
\end{equation}
It covers the group $SO(p,q)_+$, the identity component of $SO(p,q)$, consisting of the orientation preserving isometries which also preserve the orientation on any maximally space-- and timelike subspace of $\R^{p,q}$. To determine the Lie algebra of $Spin(p,q)_+$, we let $e_1,\ldots,e_{p+q}$ be an orthonormal basis of $\R^{p,q}$. The Lie algebra of $\mf{spin}(p,q)=\mf{so}(p,q)$ can be identified with the linear span of the subset $\{e_k\cdot e_l\,|\,1\leq k<l\leq p+q\}$. Thus, as a vector space 
$$
\mf{so}(p,q)=\Lambda^2\R^{p,q},
$$ 
and the commutator is given by $[x,y]=x\cdot y-y\cdot x$, where $\cdot$ denotes the algebra product of $\cliff(p,q)$.

Restricting the matrix representation of $\cliff(\C^n)$ to $Spin(p,q)$ yields the two (complex) irreducible spin representations\index{spin representation} $\Delta_{\pm}$ for $n=p+q=2m$, and the irreducible spin representation\index{spin representation} $\Delta$ for $n=2m+1$. The qualifier ``spin" refers to the fact that these representations do not factorise over $SO(p,q)$ via~(\ref{covering}), and any irreducible spin representation\index{spin representation} is isomorphic to either $\Delta_{n\pm}$, $n$ even, or $\Delta_n$, $n$ odd. Any representation of $Spin(p,q)$ that does factorise is referred to as a {\em vector representation}\index{vector representation}. Since $Spin(p,q)\subset\cliff(p,q)$, it also acts on the spaces $S_{p,q\pm}$. For instance, the spin representations of $Spin(8)$ is the complexification of $S_{8,0\pm}$ which, as vector spaces, are isomorphic to $\R^8$. In this case, one says that the spin representation is of {\em real type}. A similar analysis can be carried out for arbitrary signature using~(\ref{cliffevtable}). The induced action of $\mf{so}(p,q)$ on $\Psi\in\Delta_*$ is given on decomposable elements $x\wedge y\in\Lambda^2\R^{p,q}$ by
\begin{equation}\nonumber
x\wedge y(\Psi)=\frac{1}{4}[x,y]\cdot\Psi.
\end{equation}
In particular, for an orthonormal basis we obtain $e_k\wedge e_l(\Psi)=e_k\cdot e_l\cdot\Psi/2$. 

Finally, we want to describe a {\em $Spin(n)$--equivariant, conjugate linear} operator $\mc{A}$ which will become important in the sequel, the so--called {\em charge conjugation operator}\index{charge conjugation operator}. Again, we first assume $n=2m$ and $e_1,\ldots,e_m$ to be an orthonormal basis of $U$ in $\R^{2m}=U\oplus J(U)$. For $\Psi\in\Delta$, we let 
$$
\mc{A}(\Psi)=e_1\cdot\ldots\cdot e_m\cdot\overline{\Psi},
$$ 
where complex conjugation is defined with respect to the real form $\Lambda^*U\subset\Lambda^*U^{\C}$. The operator $\mc{A}$ thus preserves chirality if $m$ is even, and reverses chirality if $m$ is odd. Moreover, it satisfies 
$$
\mc{A}(X\cdot\Psi)=(-1)^{m+1}X\cdot \mc{A}(\Psi),\quad\mc{A}^2=(-1)^{m(m+1)/2}Id.
$$ 
The odd dimensional case simply follows by restricting $\mc{A}$ to $\cliff(\C^{2m})^{ev}$ and using the isomorphism $\gamma$ defined above. The $Spin(n)$--invariant hermitian inner product $q$ on $\Delta_n$ gives rise to the $Spin(n)$--invariant bilinear form, still written $\mc{A}$ by abuse of notation,
$$
\mc{A}(\Psi,\Phi)=q(\mc{A}(\Psi),\Phi).
$$ 
Note that $\mc{A}(\Psi,\Phi)=(-1)^{m(m+1)}\mc{A}(\Phi,\Psi)$. In particular, this shows $\Delta_n$ to be self--contragredient, i.e. its equivalent to $\Delta^*_n$ as a $Spin(n)$--representation. Furthermore, we obtain a $Spin(n)$--equivariant embedding into the $Spin(n)$--module of exterior forms defined by 
$$
[\cdot\,,\cdot]:\Delta_n\otimes\Delta_n\hookrightarrow \Lambda^*\C^n,\quad [\Psi\otimes\Phi](x_1,\ldots,x_r)=\mc{A}(\Psi,x_1\cdot\ldots\cdot x_r\cdot\Phi).
$$ 
This operation is known as {\em fierzing}\index{fierzing} in the physics' literature. Note that the map is onto for $n$ even. We will consider an example in the next paragraph.

\paragraph{Spin structures\index{spin structure}.}
Before this, we want to define spinors fields over a manifold, in analogy with vector fields. This requires the existence of a {\em spin structure\index{spin structure}}, that is, a collection of transition functions $\{\widetilde{s}_{ab}:U_a\cap U_b\to Spin(p,q)\}$ satisfying the cocycle condition~(\ref{cocycle}). On one hand side, this yields a pseudo--Riemannian vector bundle $\mb{V}^{p,q}$ of signature $(p,q)$ associated with the $SO(p,q)$--structure $s_{ab}=\pi_0\circ\widetilde{s}_{ab}$. On the other hand, we can associate the irreducible spin representations $S_{p,q(\pm)}$ or $\Delta_{(\pm)}$ via $\{\widetilde{s}_{ab}\}$ to obtain the spinor bundles 
$$
\mb{S}(\mb{V}^{p,q})_{(\pm)}=\coprod U_a\times S_{p,q(\pm)},\quad\Delta(\mb{V}^{p,q})_{(\pm)}=\coprod U_a\times\Delta_{(\pm)}.
$$
A ($\mb{V}^{p,q}$--) {\em spinor field} $\Psi$ is a section of $\mb{S}(\mb{V}^{p,q})_{(\pm)}$ or  $\Delta(\mb{V}^{p,q})_{(\pm)}$, that is, a collection of maps $\{\Psi_a:U_a\to S_{p,q(\pm)}\}$ or $\{\Psi_a:U_a\to\Delta_{(\pm)}\}$ such that $\Psi_a=\widetilde{s}_{ab}(\Psi_b)$. Conversely, assume to be given an $SO(p,q)$--structure on $M^{p+q}$. Of course, we can lift the transition functions locally to $Spin(p,q)$, but one would expect to meet topological obstructions for doing so in such a way such that the cocycle condition holds. Indeed, we have the

\begin{proposition}
Let $\mb{V}^{p,q}$ be a vector bundle associated with an $SO(p,q)_+$--valued collection of transition functions $\{s_{ab}\}$. Then there exists a spin structure\index{spin structure} which covers $\{s_{ab}\}$ if and only if the second Stiefel--Whitney class of $w_2(\mb{V}^{p,q})$ vanishes.
\end{proposition}

An important example of this is the spin structure\index{spin structure} of an oriented Riemannian manifold. In case it exists, the manifold is said to be {\em spinnable}. By the above, this is equivalent to requiring $w_2(M)=0$. It is important to realise that there might be several ways of gluing the local lifts of the transitions functions $s_{ab}$ together. Modulo equivalence, the spin structures\index{spin structure} which cover a given $SO(p,q)_+$--structure stand actually in bijection with $H^1(M,\Z_2)$ . However, in the situations we will encounter in the sequel, the spin structure\index{spin structure} will be induced by a $G$--structure $\{s_{ab}\}$ for which the inclusion $G\subset SO(p,q)$ lifts to $Spin(p,q)$, in which case a canonic spin structure\index{spin structure} is provided by $\{s_{ab}\}$ itself.

\begin{example}
Let $(M^6,g)$ be a spinnable Riemannian manifold with a given spin structure\index{spin structure}. The existence of a chiral unit spinor $\Psi$ is equivalent to an $SU(3)$\index{$SU(3)$}--structure. To see this, let us first look at the algebraic picture. The group $Spin(6)$ is isomorphic with $SU(4)$ and under this identification, the chiral spin representations $\Delta_{\pm}$ become the standard complex vector representation $\C^4$ and its complex conjugate $\overline{\C^4}$, in accordance with what we have said earlier (cf.~(\ref{cliffevtable})). A chiral (say positive) unit spinor is then an element in $S^7\subset \C^4$. But $S^7$ is acted on transitively by $SU(4)$, and in fact, isomorphic with $SU(4)/SU(3)$\index{$SU(3)$}. Hence a unit spinor on $M^6$ can be regarded as a section of the sphere bundle $\mb{S}^7=\mb{SU}(4)/\mb{SU}(3)$ with fibre $SU(4)/SU(3)$\index{$SU(3)$}, thereby inducing a reduction to $SU(3)$\index{$SU(3)$}. The inclusion $SU(3)\hookrightarrow Spin(6)$\index{$SU(3)$} induced by the choice of an $SU(3)$\index{$SU(3)$}--structure projects to an inclusion $SU(3)\hookrightarrow SO(6)$\index{$SU(3)$} via $\pi_0$~(\ref{covering}). This reduction form $SO(6)$ to $SU(3)$\index{$SU(3)$} can be understood in terms of a reduction from $SO(6)$ to $U(3)$, which yields an almost complex structure $J$ acting as an isometry for $g$, and a reduction from $U(3)$ to $SU(3)$\index{$SU(3)$}, which is equivalent to the choice of a $(3,0)$--form $\Omega=\psi_++i\psi_-$ of constant length. In particular, $\omega(X,Y)=g(JX,Y)$ defines a non--degenerate 2--form. Conversely, the embedding $SU(3)\hookrightarrow SO(6)$\index{$SU(3)$} admits a canonic lift to $Spin(6)$ since $SU(3)$\index{$SU(3)$} is simply--connected. Therefore, a reduction from $SO(6)$ to $SU(3)$\index{$SU(3)$} yields a canonic spin structure\index{spin structure} with a preferred unit spinor $\Psi$. Both points of view are intertwined by the identification $\Delta_6\otimes\Delta_6\cong\Lambda^*\C^6$ given above, namely 
$$
[\overline{\Psi}\otimes\Psi]=e^{-i\omega},\quad[\Psi\otimes\Psi]=\Omega.
$$
\end{example}

\begin{remark}
The obstruction to the existence of such a unit spinor is the Euler class $\chi\big(\Delta_6(TM)\big)$ in $H^8(M^6,\Z)$, since $\Delta_6(TM)$ is of real rank $8$. Hence the obstruction vanishes trivially, so there always exists an $SU(3)$\index{$SU(3)$}--structure if $(M^6,g)$ is spinnable.
\end{remark}

In the previous section we introduced (linear) connections, that is, a first order differential operator which act on sections of vector bundles associated with some $GL(m)$--structure. Locally, they are essentially defined by a $\mf{gl}(m)$--valued $1$--form $\omega_a$. Recall that for a Riemannian manifold, there was the canonic Levi--Civita connection whose defining $1$--form $\omega_a$ took values in $\mf{so}(p,q)$. Consequently, the Levi--Civita connection acts on the spin representation, and we can covariantly derive spinors, too.

\begin{example}
Let us take up again the previous example. Since the fierzing map $[\cdot,\cdot]$ is $Spin(6)$--equivariant, it commutes with $\nabla^{LC}$. Therefore, if $\nabla^{LC}\Psi=0$, then $\nabla^{LC}\omega=0$ and $\nabla^{LC}\Omega=0$. As a result, we have a geometrical reduction to $SU(3)$\index{$SU(3)$}, or equivalently, the {\em holonomy} of the metric reduces to $SU(3)$\index{$SU(3)$}.
\end{example}

The Levi--Civita connection on the spinor bundle $\mb{S}_{p,q(\pm)}$ gives also rise to several differential operators. One we will encounter frequently is the {\em Dirac operator\index{Dirac operator}} $\mb{D}$ of the spin structure\index{spin structure}. Using the metric to identify $TM$ with its dual $T^*M$, we can define $\mb{D}$ as
$$
\mb{D}:\Gamma(\mb{S}_{p,q(\pm)})\stackrel{\nabla^{LC}}{\longrightarrow}\Gamma(T^*M\otimes\mb{S}_{p,q(\pm)})\stackrel{g}{\longrightarrow}\Gamma\big(TM\otimes\mb{S}_{p,q(\pm)}\big)\stackrel{\cdot}{\longrightarrow}\Gamma(\mb{S}_{p,q(\mp)}),
$$
where $\cdot$ denotes again Clifford multiplication\index{Clifford multiplication}.
%
%
\subsection{Compactification in supergravity}
\label{sugracom}

For the moment being, we know five consistent supersymmetric string theories, namely a so--called {\em type I\index{type I}} string theory, two {\em heterotic}\index{heterotic} string theories and two {\em type II}\index{type II} string theories. They are all defined over a ten--dimensional space--time $M^{1,9}$, but apart from this, their mathematical formulation has little else in common at first glance. However, they are all supposed to give rise to the same observable 4--dimensional physics as reflected in the existence of various dualities between these string theories. The low energy limit of the string theory gives rise to the corresponding supergravity theory, and our aim is to give a $G$--structure interpretation of type II\index{type II} supergravity. For this, it is instructive to understand how $G$--structures emerge in heterotic\index{heterotic} supergravity first.

\paragraph{Heterotic\index{heterotic} supergravity.}
A good (yet not exhaustive) list of references for the material of this section is provided by~\cite{ccdlmz03},~\cite{gmw04},~\cite{chsa02} and~\cite{iviv05}. 

The field content on $M^{1,9}$ of the two heterotic\index{heterotic} supergravity theories consists of
\begin{itemize}
\item a {\em space--time metric} $g$.
\item a {\em dilaton field}\index{dilaton field} $\phi\in C^{\infty}(M)$.
\item an {\em $H$--flux} $H\in\Omega^3(M)$. 
\item a {\em gauge field} $F$ in $\Omega^2(M,\mf{e}_8\times\mf{e}_8)$ or $\Omega^2\big(M,\mf{so}(32)/\Z_2\big)$, depending on which type of heterotic\index{heterotic} theory we consider.
\item a {\em supersymmetry parameter\index{supersymmetry parameters}} $\Psi$, a chiral spinor of unit norm. 
\end{itemize}
Moreover, this datum is supposed to satisfy the following field equations:
\begin{itemize} 
\item the {\em modified Bianchi identity} $dH=2\alpha'Tr(F\wedge F)$ ($\alpha'$ being a universal constant, the {\em string torsion}).
\item the {\em gravitino equation\index{gravitino equation}} $\nabla^{LC}_X\Psi+\frac{1}{4}(X\llcorner H)\cdot\Psi=0$.
\item the {\em dilatino equation}\index{dilatino equation} $(d\phi+\frac{1}{2}H)\cdot\Psi=0$.
\item the {\em gaugino equation} $F\cdot\Psi=0$.
\end{itemize}
In order to find a solution, one usually makes a {\em compactification} ansatz, that is, one considers a space--time of the form $(M^{1,9},g^{1,9})=(\R^{1,3},g_0)\times(M^6,g)$, where $(\R^{1,3},g_0)$ is flat Minkowski space and $(M^6,g)$ a $6$--dimensional spinnable Riemannian manifold. This is not only a convenient mathematical ansatz, it also reflects the empiric fact that the phenomenologically tangible world is confined to three spatial dimensions plus time. Then, one tries to solve the above equations with fields living on $M^6$, trivially extended to the entire space--time. Let us assume $F=0$ for further simplification, so we are looking for a set of datum $(g,\phi,H,\Psi)$ (with $\Psi$ a section of, say, $\Delta_+(TM)$), that satisfies the gravitino\index{gravitino equation} and the dilatino equation\index{dilatino equation} on $M$. As the previous example shows, the (spinnable) metric induces a reduction to $SO(6)$, and the unit spinor $\Psi$ yields a further reduction from $Spin(6)$ to $SU(3)$\index{$SU(3)$}. The field equations can now be interpreted as a differential constraint on the $SU(3)$\index{$SU(3)$}--structure. For this, we first have to analyse the action of the Levi--Civita connection which we think of as a collection of local differential operators ${d+\omega_a}$, with $\omega_a$ a $1$--form taking values in $\mf{so}(6)=\mf{su}(4)$. On the other hand, the spinor $\Psi$ is locally given by a {\em constant} map $\Psi_a:U_a\to\Delta_+$, whence 
$$
(\nabla_X^{LC}\Psi)_a=\omega_a(X_a)(\Psi_a).
$$
Now $\mf{su}(4)=\mf{su}(3)\oplus\mf{su}(3)^{\perp}$, and since $\mf{su}(3)$ acts trivially on $\Psi$, the action of the connection is encapsulated in the tensor  
$$
\mc{T}=proj_{\mf{su}(3)^{\perp}}(\omega_a)\in\Gamma(\Lambda^1T^*M\otimes\mf{su}(3)^{\perp})\subset\Gamma(\Lambda^1T^*M\otimes\Lambda^2T^*M),
$$
the so--called {\em intrinsic torsion}. It is a first order differential geometric invariant and measures the failure of the local sections $s_a$ of $TM$ which define the $SU(3)$\index{$SU(3)$}--valued transition functions, to be induced by a coordinate system~\cite{be60}. The gravitino equation\index{gravitino equation} therefore states that we are looking for an $SU(3)$\index{$SU(3)$}--structure whose intrinsic torsion is skew--symmetric (algebraic constraint) and closed (topological constraint). The algebraic condition can be analysed by using representation theory. Skew--symmetry of the intrinsic torsion means that at any point, $\mc{T}_p$ lies in the image of the equivariant embedding $\Lambda^3T^*_pM\hookrightarrow T^*_pM\otimes\Lambda^2T^*_pM$ followed by projection onto $T^*_pM\otimes\mf{su}(3)^{\perp}$. Decomposing this space and $\Lambda^3T^*_pM$ into $SU(3)$\index{$SU(3)$}--irreducibles yields
$$
\Lambda^3T_p^*M=2\R\oplus\llbracket\C^3\rrbracket\oplus\llbracket \odot^{2,0}\rrbracket,\quad T^*M\otimes\mf{su}(3)^{\perp}=2\R\oplus2\llbracket \C^3\rrbracket\oplus2\mf{su}(3)\oplus\llbracket \odot^{2,0}\rrbracket,
$$
where $\R$ denotes the trivial representation and $\mf{su}(3)$ the adjoint representation of $SU(3)$\index{$SU(3)$} both of which are real and occur with multiplicity $2$, $\llbracket\C^3\rrbracket$ the real representation obtained from the standard representation of $SU(3)$\index{$SU(3)$} on $\C^3$ by forgetting the complex structure, and $\llbracket \odot^{2,0}\rrbracket$ the real representation obtained by forgetting the complex structure on the $SU(3)$\index{$SU(3)$}--representation of complex symmetric $2$--tensors over $\C^3$. Because of equivariance, Schur's Lemma\index{Schur's Lemma} implies that the only modules of $T^*_pM\otimes\mf{su}(3)^{\perp}$ which can be hit by a $3$--form are $\R$, $\llbracket\C^3\rrbracket$ and $\llbracket\odot^{2,0}\rrbracket$, as otherwise, there would be non--trivial equivariant maps between non--equivalent $SU(3)$\index{$SU(3)$}--representation spaces. In particular, we see at once that $\mc{T}_p$ is not allowed to take values in $\mf{su}(3)\oplus\mf{su}(3)$. The dilatino equation\index{dilatino equation} can be discussed in a similar vein and forces a certain component of the intrinsic torsion to be exact. The decisive advantage of using an $SU(3)$\index{$SU(3)$}--structure ansatz here is that the intrinsic torsion can be computed out of the differentials of the $SU(3)$\index{$SU(3)$}--invariant forms $\psi_{\pm}$ and $\omega$ via the fierzing map $[\cdot\,,\cdot]$, as this map is $Spin(6)$--equivariant. We considered a special case of this before, when we observed that $\nabla^{LC}\Psi=0$ is equivalent to $\nabla\omega=\nabla\Omega=0$, which turns out to hold precisely if $d\omega=d\Omega=0$. More generally, as $\mc{T}$ is an invariant of the $SU(3)$\index{$SU(3)$}--structure, the covariant derivative of the $SU(3)$\index{$SU(3)$}--{\em invariant} differential forms $\omega$ and $\psi_{\pm}$ are
$$
\nabla^{LC}_X\omega=\mc{T}(X)(\omega)\quad\nabla^{LC}_X\psi_{\pm}=\mc{T}(X)(\psi_{\pm}).
$$ 
Now the exterior differential operator on forms is the skew--symmetrisation of $\nabla$, that is, for a local basis we have $d=\sum s_k\wedge\nabla^{LC}_{s_k}$. It follows that the exterior differentials $d\omega$ and $d\psi_{\pm}$ are determined by the intrinsic torsion. It turns out that, although the exterior derivative contains considerably less information than the covariant derivative, the intrinsic torsion is completely captured by $d\omega$ and $d\psi_{\pm}$. This fact is a special feature of $SU(3)$\index{$SU(3)$}--structures and by no means true in general, though similar results hold for other structure groups in low dimensions. Concretely, this means that we translated the spinor field equations into integrability conditions on $d\omega$ and $d\psi_{\pm}$. In practice, the differentials are far easier to compute than the spinor derivative and in this way, explicit solutions to heterotic\index{heterotic} string theory could be found~\cite{ccdlmz03}.

\paragraph{Type II\index{type II} supergravity.}\label{typeIIsugra}

For a brief, mathematically flavoured review of type II\index{type II} supergravity see~\cite{jewi05b} and the references quoted there for details.

Type II\index{type II} supergravity requires two supersymmetry parameters\index{supersymmetry parameters} $\Psi_L$ and $\Psi_R$ (whence ``type II\index{type II}''), which are unit spinors of equal (type IIB\index{type IIB}) or opposite (type IIA\index{type IIA}) chirality. The fields on $M^{1,9}$ come in two flavours; they are either NS-NS\index{NS-NS} (NS=Neveu-Schwarz\index{Neveu--Schwarz}) or R-R\index{R-R} (R=Ramond). To the former class belong
\begin{itemize}
\item the space--time metric $g$.
\item the {\em $B$--fields}, a collection of locally defined $2$--forms $\{B_a\in\Omega^2(U_a)\}$, whose differentials glue to the closed 3--form $H$--{\em flux} $H=dB_a$. Moreover, quantisation arguments require $H$ to be integral.
\item the dilaton field\index{dilaton field} $\phi\in C^{\infty}(M)$.
\end{itemize}
The R-R sector\index{R-R sector} consists of a closed differential form $F$ of either even (type IIA\index{type IIA}) or odd degree (type IIB\index{type IIB}). The homogeneous components of $F$ are referred to as {\em Ramond--Ramond fields}\index{Ramond--Ramond field}. Only half of these are physical in the sense that they represent independent degrees of freedom, so we have a duality relation 
$$
F^p=(-1)^{p(p+1)/2}\star F^{10-p}.
$$
Locally, $F=dC_a$ for a collection of forms $\{C_a\in\Omega^{ev,od}(U_a)\}$ of suitable parity, the {\em Ramond--Ramond potentials}. The field equations are encapsulated in the so--called {\em democratic formulation}\index{democratic formulation} of Bergshoeff et al.~\cite{bkop01}, which for instance compactified to six dimensions yields for type IIA\index{type IIA}\footnote{Similar equations hold for type IIB\index{type IIB}.} the gravitino equation\index{gravitino equation}
\begin{equation}\label{mastereqgrav}
\begin{array}{lcl}
\nabla_X\Psi_L - \frac{1}{4}(X\llcorner H)\cdot \Psi_L
-e^{\phi}F\cdot X\cdot \mc{A}(\Psi_R) & = & 0\\[5pt]
\nabla_X\Psi_R +\frac{1}{4}(X\llcorner H)\cdot\Psi_R+e^{\phi}F\cdot X\cdot\mc{A}(\Psi_L) & = & 0,
\end{array}
\end{equation}
and the dilatino equation\index{dilatino equation},
\begin{equation}\label{mastereqdil}
\begin{array}{c}
({\rm D} - d\phi - \frac{1}{4} H)\cdot\Psi_L=0,\quad
({\rm D} - d\phi + \frac{1}{4} H)\cdot\Psi_R=0,
\end{array}
\end{equation}
involving the Dirac operator\index{Dirac operator} of the spin structure\index{spin structure}.

One would like to discuss type II\index{type II} supergravity along the lines of heterotic\index{heterotic} supergravity. Are classical $G$--structures of any use here? No. Let us see where the problem occurs. Mimicking the approach of the previous section, we are looking for a $G$--structure where $G$ stabilises two spinors $\Psi_L$ and $\Psi_R$. In $6$ dimensions, that stabiliser inside $Spin(6)$ is $SU(2)$. Since $\Psi_L$ and $\Psi_R$ are $SU(2)$--invariant, so is their angle $q(\Psi_L,\Psi_R)$. However, the physical model allows for totally independent unit spinors, so we rather need two independent $SU(3)$\index{$SU(3)$}--structures. In general they do not intersect in a well--defined substructure. For instance, the spinors might coincide at some points where both stabilisers would intersect in $SU(3)$\index{$SU(3)$} while outside the coincidence set they would pointwise intersect in $SU(2)$ (though this does still not imply a global reduction to $SU(2)$). The geometric configuration of type II\index{type II} supergravity requires therefore a vector bundle other than the tangent bundle. This is where generalised geometry enters the scene.
%
%
%
%
%
\section{The linear algebra of $Spin(n,n)$}\label{linalg}
As we saw in the previous section, the possible reductions from a group $\mb{G}$ to a subgroup $G$ are parametrised by the orbits $\mb{G}/G$. It is therefore essential to study first the coset space $\mb{G}/G$ which parametrises the $G$--structures on the representation space of $\mb{G}$. For instance, a Euclidean vector space $(\R^n,g)$ can be thought of as the result of reducing the structure group $\mb{G}=GL(n)$ of $\R^n$ to $G=O(n,g)$, that is, the choice of a Euclidean metric $g$ singles out an embedding of $O(n)$ into $GL(n)$, whose image is $O(n,g)$. The starting point for generalised geometry are the groups $\mb{G}=SO(n,n)$ and $Spin(n,n)$.
%
%
\subsection{The group $Spin(n,n)$}\label{sonn}
In Section~\ref{spinstruc} we constructed a representation space of $\cliff(\C^n)$ as the exterior algebra over some subspace of $\C^n$. Similarly, we can construct a representation space of $\cliff(n,n)$, which involves the choice of a splitting $\R^{n,n}=W\oplus W'$ into two maximally isotropic subspaces. In fact, this choice comes effectively down to an isometry between $(\R^{n,n},g)$ and $\big(W\oplus W^*,(\cdot,\cdot)\big)$ with contraction as inner product, i.e. $(w,\xi)=\xi(w)/2$ for $w\in W$ and $\xi\in W^*$: The map $w'\mapsto g(w',\cdot)/2\in W^*$ is injective since $W$ is isotropic\footnote{The factor $1/2$ is introduced for computational purposes and has no geometrical meaning. Note also that we use conventions slightly different from~\cite{hi03} and~\cite{wi06} which results in different signs and scaling factors.}. Moreover, the choice of a preferred isotropic subspace gives $\R^{n,n}$ a preferred orientation. Namely, $W$ defines a subgroup $GL(W)\subset O(n,n)$, as for $A\in GL(W)$, $$
(Aw,A^*\xi)=\frac{1}{2}A^*\xi(Aw)=\frac{1}{2}\xi(A^{-1}Aw)=(w,\xi).
$$ 
The two connected components of $GL(W)$ single out the connected components of $O(n,n)$ which make up the group $SO(n,n)$. If we give $W$ itself an orientation, so that the structure group of $W$ is reduced to the identity component $GL(W)_+$, this argument also shows the structure group of $W\oplus W^*$ to reduce to its identity component $SO(n,n)_+$.

As a vector space, the spin representation is $P_{n,n}=\Lambda^*W^*$ and an element $x\oplus \xi\in W\oplus W^*$ acts on $\rho\in P_{n,n}$ by 
$$
(x\oplus \xi)\bullet\rho=-x\llcorner\rho+\xi\wedge\rho.
$$ 
Indeed, one easily checks that $(x\oplus \xi)^2=-(x,\xi)Id$, so that in virtue of the universal property, the map 
$$
x\oplus \xi\in W\oplus W^*\mapsto -x\llcorner+\xi\wedge\in\End(P_{n,n})
$$
extends to an isomorphism $\cliff(n,n)\cong\End(P)$, in accordance with~(\ref{clifftable}) (where from now on we drop the subscript $n,n$ to ease notation). Restricting this action to $Spin(n,n)$ yields the irreducible spin representations $S_{\pm}=\Lambda^{ev,od}W^*$. The inclusion $GL(W)\hookrightarrow SO(n,n)$ can be lifted to $Spin(n,n)$, albeit in a non--canonic way. In practice, we always assume to have chosen an orientation on $W$, as $GL(W)_+$ naturally lifts to $Spin(n,n)$. Restricted to this lift, the spin representation becomes the exterior algebra tensored with the square root of the line bundle spanned by $n$--vectors,
\begin{equation}\label{spinGL}
S_{\pm}\cong\Lambda^{ev,od}W^*\otimes\sqrt{\Lambda^nW}\mbox{ as a }GL(W)_+\subset Spin(n,n)-\mbox{space},
\end{equation}
a fact which will be important to bear in mind later on.

There is also a $Spin(n,n)$--invariant bilinear form on $S_{\pm}$. Let $\wedge$ be the anti--automorphism defined on algebra elements $a^p$ of degree $p$ by $\widehat{a^p}=(-1)^{p(p+1)/2}$. On $S_{\pm}$, we define 
$$
\langle\rho,\tau\rangle=[\rho\wedge\widehat{\tau}]^n\in \Lambda^nW^*,
$$
where $[\cdot]^n$ denotes projection on forms of degree $n$. After choosing a non--zero volume form on $W$, i.e. a trivialisation of $\Lambda^nW^*$, this form takes values in the reals. It is symmetric for $n\equiv 0,3\!\mod4$ and skew for $n\equiv 1,2\!\mod4$, i.e. 
$$
\langle\rho,\tau\rangle=(-1)^{n(n+1)/2}\langle\tau,\rho\rangle.
$$ 
Moreover, $S_+$ and $S_-$ are non--degenerate and orthogonal if $n$ is even and totally isotropic if $n$ is odd. 

A particularly important subset of transformations in $Spin(n,n)$ is given by the so--called {\em $B$--field transformations}\index{$B$--field transformation}. As a $GL(W)$--space, 
$$
\mf{so}(n,n)\cong\Lambda^2(W\oplus W^*)=\Lambda^2W\oplus W\otimes W^*\oplus\Lambda^2W^*
$$
which shows that any 2--form $B=\sum B_{kl}w^k\wedge w^l$ over $W$ acts through exponentiation  as an element of both $Spin(n,n)_+$ and $SO(n,n)_+$. Concretely, $B$ becomes a skew--symmetric operator on $W\oplus W^*$ via the embedding
$$
\zeta\wedge\eta(x\oplus\xi)=(\zeta,X\oplus\xi)\xi-(\eta,X\oplus\xi)\zeta=\frac{1}{2}X\llcorner(\zeta\wedge\eta).
$$
Then
$$
e^B_{SO(n,n)}(X\oplus\xi)=\left(\begin{array}{cc} 1\,\, & 0\\ B/2\,\, & 1\end{array}\right)\left(\begin{array}{c} X\\ \xi\end{array}\right),
$$
where the matrix is taken with respect to the splitting $W\oplus W^*$. On the other hand, $B\in\Lambda^2(W\oplus W^*)$ sits naturally inside $\cliff(W\oplus W^*)$. Again by standard representation theory, the exponential of $B$ acts on a spinor $\rho$ via
$$
e^B\bullet\rho=(1+B+\frac{1}{2}B\bullet B+\ldots)\bullet\rho=(1+B+\frac{1}{2}B\wedge B+\ldots)\wedge\rho=e^B\wedge\rho.
$$
Note that the differential of the covering map $\pi_0:Spin(n,n)\to SO(n,n)$ links the exponentials via
\begin{equation}\label{diffB}
\pi_0(e^B_{Spin(n,n)})=e^{\pi_{0*}(B)}_{SO(n,n)}=e^{2B}_{SO(n,n)}.
\end{equation}
%
%
\subsection{Special orbits}
\label{specorb}
In the introduction we encountered generalised Riemannian and generalised $SU(3)$\index{$SU(3)$}--structures. Next we describe these and various other generalised structures in terms of structure groups, that is, we wish to see these structures as being defined by sections of a fibre bundle whose fibre is a coset space of the form $SO(n,n)_+/G$ or $Spin(n,n)_+/G$ (cf. Section~\ref{gstructures}). 

To begin with, consider a straight generalised Riemannian structure on $W\oplus W^*$ induced by some Riemannian metric $g$ on $W$\footnote{More accurately, one should speak of a straight generalised Euclidean structure being induced by a Euclidean structure on the oriented vector space $(W,g)$ as we are doing linear algebra for the moment, but we do not wish to overload the terminology here.}. It is defined by the linear endomorphism
$$
\mc{G}_0=\left(\begin{array}{cc} 0 & g^{-1}\\ g & 0\end{array}\right)
$$
whose matrix is taken with respect to the decomposition $W\oplus W^*$. This squares to the identity, and the $\pm1$--eigenspaces $D^{\pm}$ are given by
$$
D^{\pm}=\{X\oplus\pm g(X,\cdot)\,|\,X\in W\}.
$$
Moreover, the restriction of $(\cdot\,,\cdot)$ to $D^{\pm}$ induces a positive and negative definite inner product $g_+$ and $g_-$ on $D^+$ and $D^-$. This orthogonal decomposition
$$
\big(W\oplus W^*,(\cdot\,,\cdot)\big)=\big(D^+\oplus D^-,g_+\oplus g_-\big)
$$
is preserved by the group 
\begin{eqnarray}
SO(n,0)\times SO(0,n)\!\! &\!\! \cong \!\!&\!\! \{\left(\begin{array}{cc} \!\!A_+ & \!\!0\\\!\! 0 & \!\!A_-\end{array}\right)| A_{\pm}\in\End(D^{\pm}),\,A^*_{\pm}g_{\pm}=g_{\pm},\det A_{\pm}=1\},\nonumber\\
\!\!&\!\! =\!\! & \!\!SO(D^+)\times SO(D^-)\label{straightembed}
\end{eqnarray}
leading to the following

\begin{definition}
A {\em generalised Riemannian structure} on $W\oplus W^*$ is the choice of an embedding $SO(n,0)\times SO(0,n)\hookrightarrow SO(n,n)_+$.
\end{definition}

For the straight case, this subgroup is given by the embedding~(\ref{straightembed}) and corresponds to the endomorphisms $\mc{G}_0$. The $B$--field transform $\mc{G}_B$ then corresponds to conjugation with $\exp(2B)$, i.e. to $\exp(2B)\big(SO(D^+)\times SO(D^-)\big)\exp(-2B)$. Actually, any generalised Riemannian metric is the $B$--field transform of a straight generalised Riemannian structure (a fact wrong for other generalised structures). For this we note that by definition any subgroup of the form $SO(n,0)\times SO(0,n)$ arises as the stabiliser of an orthogonal decomposition $\big(V^+\oplus V^-,g_+\oplus g_-\big)$ of $\big(W\oplus W^*,(\cdot\,,\cdot)\big)$. The definite spaces $V^{\pm}$ intersect the isotropic spaces $W$ and $W^*$ trivially, so they can be written as the graph of an isomorphism $P^{\pm}:W\to W^*$. Dualising $P^{\pm}$ yields an element in $W^*\otimes W^*$ with $P^-=-P^{+tr}$. The symmetric part defines a Riemannian metric $g=(P^++P^{+tr})/2$, while the skew--symmetrisation yields a 2--form $B=(P^+-P^{+tr})/2$ on $W$. The endomorphism $\mc{G}$ defined by $\mc{G}_{|V^{\pm}}=\pm Id_{V^{\pm}}$ then coincides with $\mc{G}_B$, whence the

\begin{proposition} 
A {\em generalised Riemannian metric} on $W\oplus W^*$ is characterised by either of the following, equivalent statements:

\noindent{\rm (i)} The structure group reduces from $SO(n,n)_+$ to $SO(n)\times SO(n)$. 

\noindent{\rm (ii)} The choice of a pair $(g,B)$, consisting of a positive definite inner product $g\in\odot^2W^*$ and a $2$--form $B\in\Lambda^2W^*$.

\noindent{\rm (iii)} The choice of a positive definite, oriented subspace $V^+\subset W\oplus W^*$ which is of maximal rank, i.e. ${\rm rk}\,V^+=\dim W$. 
\end{proposition}

The involution $\mc{G}$ corresponding to the generalised metric preserves the metric, but reverses the orientation if $n$ is odd. At any rate, it defines an element in $O(n,n)$ and as such, it can be lifted to an element $\gtilde$ in $Pin(n,n)$. We analyse the action of $\gtilde$ on spinors next. First consider the case of $\mc{G}$ being induced by $D^+\oplus D^-$ with oriented orthonormal basis $d_k^{\pm}=e_k\oplus\pm g(e_k,\cdot)$ ($e_j$, $j=1,\ldots,n$ being an orthonormal basis of $(W,g)$). Then $\mc{G}$ is the composition of reflections $R_{d^-_k}$ along $d_k^-$,
$$
\mc{G}=R_{d^-_1}\circ\ldots\circ R_{d^-_n}.
$$ 
Therefore, $\widetilde{\mc{G}}$ acts via Clifford multiplication\index{Clifford multiplication} as the Riemannian volume form $vol_{D^-}=d^-_1\wedge\ldots\wedge d^-_n$ of $D^-$. Next we express the $\wedge$--product and the Hodge $\star$--operator of $\Lambda^*W^*$ in terms of the Clifford algebra product in $\cliff(W,g)$ via the natural isomorphism $\mf{J}:\cliff(W,g)\to\Lambda^*W^*$. For any $x\in W$ and $a\in \cliff(W,g)$, 
$$
\mf{J}(x\cdot a)=-x\llcorner\mf{J}(a)+x\wedge \mf{J}(a),\quad \mf{J}(a\cdot x)=(-1)^{\deg(a)}\big(x\llcorner \mf{J}(a)+x\wedge \mf{J}(a)\big),
$$
with $\llcorner$ the adjoint of $\wedge$ with respect to $g$ (i.e. $\llcorner$ is metric contraction), and 
$$
\mf{J}(\widehat{a}\cdot vol_g)=\star_g\mf{J}(a)
$$ 
for the Riemannian volume form $vol_g$ on $W$. Moreover, $vol_g\cdot a=\widetilde{a}\cdot vol_g$ for $n$ even and $vol_g\cdot a=a\cdot vol_g$ for $n$ odd, where $\sim$ denotes the involution defined on elements of even or odd degree by $\pm id$. For $\rho\in S_{\pm}$, we thus obtain
$$
vol_{D^-}\bullet \rho=\mf{J}\big( \mf{J}^{-1}(\rho)\cdot vol_g\big).
$$ 
This implies $vol_{D^-}\bullet\rho=\star\widehat{\rho}$ if $n$ is even and $vol_{D^-}\bullet\rho=-\star\widehat{\widetilde{\rho}}$ for $n$ odd. For a non--trivial B--field, $\mc{G}$ gets conjugated by $\exp(2B)$ and thus $\widetilde{\mc{G}}$ by $\exp(B)$ (cf.~(\ref{diffB})).

\begin{proposition}
The operator $\gtilde=vol_{V^-}$ corresponding to the generalised metric $(g,B)$ acts on $S_{\pm}$ via
$$
\gtilde\bullet\rho =\left\{\begin{array}{ll} e^B\bullet\star_g(e^{-B}\bullet\rho)^{\wedge},\quad&n\mbox{ even}\\
e^B\bullet\star_g(e^{-B}\bullet\widetilde{\rho})^{\wedge},\quad&n\mbox{ odd}\end{array}\right..
$$
\end{proposition}

Up to signs, $\gtilde$ coincides with the $\Box$--operator in~\cite{wi06}. Note that 
$$
\gtilde^2=(-1)^{n(n-1)/2},\quad\langle\gtilde\bullet\rho,\tau\rangle=(-1)^{n(n+1)/2}\langle\rho,\gtilde\tau\rangle.
$$
In particular, $\widetilde{\mc{G}}$ defines a complex structure on $P$ if $n\equiv2,3\!\mod4$. 

The presence of a generalised metric also implies a very useful description of the complexification $S_{\pm}\otimes\C$ as a tensor product of the complex spin representations of $Spin(n)$. The orthogonal decomposition of $W\oplus W^*$ into $V^+\oplus V^-$ implies two things: Firstly, we can lift any vector $x\in W$ to $x^{\pm}=x\oplus P^{\pm}x\in V^{\pm}$. Secondly, $\cliff(W\oplus W^*)$ is isomorphic with the twisted tensor product $\cliff(V^+)\htimes\cliff(V^-)$\footnote{The twisted tensor product $\htimes$ of two graded algebras $A$ and $B$ is defined on elements of pure degree as $a\htimes b\cdot a'\htimes b'=(-1)^{\deg(b)\cdot\deg(a')} a\cdot a'\htimes b\cdot b'$}. We get an isomorphism between $\cliff(W,g)\htimes\cliff(W,-g)$ and $\cliff(V^+)\htimes\cliff(V^-)$ by extending
$$
x\htimes y\in W\htimes W\mapsto x^+\bullet y^-.
$$
Next we inject the $Spin(n)\times Spin(n)$--module $\Delta_n\otimes\Delta_n$ into $\Lambda^*W^*\otimes\C$ via the fierzing map\index{fierzing} $[\cdot\,,\cdot]$ introduced in~\ref{spinstruc}. As noted before, this is an isomorphism for $n$ even; in the odd case, we obtain an isomorphism by concatenating $[\cdot\,,\cdot]$ with projection on the even or odd forms, which we write as $[\cdot\,,\cdot]^{ev,od}$. We incorporate the $B$--field by defining $[\cdot\,,\cdot]_B=e^B\bullet[\cdot\,,\cdot]$. A vector $x\in W$ acts on $\Delta_n$ via Clifford multiplication\index{Clifford multiplication}, and on $S_{\pm}\otimes\C$ via Clifford multiplication\index{Clifford multiplication} with $x^{\pm}$. The next proposition states how these actions relate under $[\cdot\,,\cdot]$.

\begin{proposition}\label{commutation} {\rm~\cite{gmwi06},~\cite{wi06}} We have
\begin{eqnarray*}
\,[x\cdot\Psi_L\otimes\Psi_R]_B & = & (-1)^{n(n-1)/2}x^+\bullet[\Psi_L\otimes\Psi_R]_B ,\\
\,[\Psi_L\otimes y\cdot\Psi_R]_B & = & -y^-\bullet\widetilde{[\Psi_L\otimes\Psi_R]}_B.
\end{eqnarray*}
\end{proposition}

\begin{remark}
The sign twist induced by $\sim$ is a result of considering $\Delta\otimes\Delta$ as a $Spin(n,0)\times Spin(n,0)$-- rather than as a $Spin(n,0)\times Spin(0,n)$--module: In $\cliff(W,-g)$, an element $x$ of unit norm squares to $1$ instead of $-1$ which is compensated precisely by $\sim$. 
\end{remark}

By the preceding remark, we obtain as a

\begin{coro}\label{spincommut}
The map $[\cdot\,,\cdot]:\Delta_n\otimes\Delta_n\to S_{\pm}\otimes\C$ is $Spin(W,g)\times Spin(W,-g)$--equivariant, where this group acts on $S_{\pm}\otimes\C$ via the embedding $x\htimes y\mapsto x^+\bullet y^-$.
\end{coro}

Using Proposition~\ref{commutation}, we see that bi--spinors are ``self--dual" in the following sense: we have
$$
\gtilde\bullet[\Psi_L\otimes\Psi_R]_B=(-1)^m[\Psi_L\otimes vol_g\cdot\Psi_R]
$$ 
for $n=2m$ and 
$$
\gtilde\bullet[\Psi_L\otimes\Psi_R]_B=(-1)^m\widetilde{[\Psi_L\otimes vol_g\cdot\Psi_R]}
$$ 
for $n=2m+1$. By a standard result, the action of $vol_g$ on chiral spinors is given by
$$
\begin{array}{lcll}
vol_g\cdot\Psi_{\pm} & = & \pm(-1)^{m(m+1)/2}i^m\Psi_{\pm}\quad & \mbox{ if }n=2m\\
vol_g\cdot\Psi & = & (-1)^{m(m+1)/2}i^{m+1}\Psi\quad & \mbox{ if }n=2m+1,
\end{array}
$$
from which we deduce the

\begin{coro}\label{volume}
Let $\Psi_{L,R}\in\Delta$. 

\noindent{\rm (i)} If $n=2m$ and $\Psi_{L,R}$ are chiral, then for $\Psi_R\in\Delta_{\pm}$
$$
\gtilde\bullet[\Psi_L\otimes\Psi_R]_B=\pm(-1)^{m(m-1)/2}i^m[\Psi_L\otimes\Psi_R]_B.
$$
\noindent{\rm (ii)} If $n=2m+1$, then
$$
\gtilde\bullet[\Psi_L\otimes\Psi_R]_B=(-1)^{m(m-1)/2}i^{m+1}\widetilde{[\Psi_L\otimes\Psi_R]}_B.
$$
\end{coro}

Corollary~\ref{spincommut} opens the way to envisage further reductions down to subgroups in $Spin(n,0)\times Spin(0,n)$. Let $G_L,G_R\subset Spin(n)$ be the stabiliser of a collection of spinors $\{\Psi_{L,k}\}_k$ and $\{\Psi_{R,l}\}_l$. By the corollary, the stabiliser of the collection of $W\oplus W^*$--spinors 
$$
\{\rho_{kl}=[\Psi_{L,k}\otimes\Psi_{R,l}]_B\}
$$ 
inside $Spin(n,0)\times Spin(0,n)$ is isomorphic with $G_L\times G_R$\footnote{Where we think of $G_R$ as a subgroup of $Spin(0,n)$ under the canonical isomorphism $Spin(n,0)\cong Spin(0,n)$}. The set of datum $(g,B,\{\Psi_{L,k}\}_k,\{\Psi_{R,l}\}_l)$ therefore induces a reduction to $G_L\times G_R$ inside $Spin(n,0)\times Spin(0,n)$. Conversely, such a reduction induces the set of spinors, as we can project $G_L\times G_R$ down to $SO(V^+)\times SO(V^-)$. As $(V^{\pm},g_{\pm})$ are isometric to $(W,\pm g)$, we can pull back the $G_{L,R}$--structure on $V^{\pm}$ to $(W,\pm g)$, where the lift to $Spin(W,\pm g)$ yields again the groups we started with. Here are some examples.

\paragraph{Generalised $SU(m)$--structures\index{generalised $SU(m)$--structure}~\cite{gmwi06},~\cite{jewi05a}.}
These structures are defined by the choice of an embedding $SU(m)\times SU(m)\hookrightarrow Spin(2m,0)\times Spin(0,2m)$. The group $SU(m)\subset Spin(2m)$ stabilises two spinors $\Psi$ and $\mc{A}(\Psi)$, which are orthogonal to each other. They are thus of equal chirality for $m$ even and of opposite chirality for $m$ odd (cf. Section~\ref{spinstruc}). In the former case, we have two conjugacy classes inside $Spin(2m)$, stabilising a pair of spinors of positive or negative chirality respectively. We will always assume that the embedded $SU(m)$\index{$SU(m)$} belongs to the conjugacy class fixing a pair of positive spinors for $m$ even. The set of $W\oplus W^*$--spinors
$$
([\mc{A}(\Psi_L)\otimes\Psi_R]_B,[\Psi_L\otimes\Psi_R]_B,[\Psi_L\otimes\mc{A}(\Psi_R)]_B,[\mc{A}(\Psi_L)\otimes\mc{A}(\Psi_R)]_B)
$$ 
induces a reduction to $SU(m)_L\times SU(m)_R$\index{$SU(m)$}. Furthermore, 
$$
[\Psi_L\otimes\mc{A}(\Psi_R)]=(-1)^{m(m+1)/2}\overline{[\mc{A}(\Psi_L)\otimes\Psi_R]},
$$ 
so the reduction is already characterised by the pair 
$$
(\rho_0,\rho_1)=([\mc{A}(\Psi_L)\otimes\Psi_R]_B,[\Psi_L\otimes\Psi_R]_B).
$$
While $[\mc{A}(\Psi_L)\otimes\Psi_R]$ is always of even parity, $\widetilde{[\Psi_L\otimes\Psi_R]}=(-1)^m[\Psi_L\otimes\Psi_R]$. Moreover, $\rho_{0,1}$ are self--dual in the sense that
\begin{equation}\label{selfdualgensum}
\gtilde\bullet\rho_{0,1}=(-1)^{m(m-1)/2}i^m\rho_{0,1}.
\end{equation} 
The name of generalised $SU(m)$--structure\index{generalised $SU(m)$--structure} is justified by the fact that this class comprises the $B$--field transformations\index{$B$--field transformation} of ``classical'' $SU(m)$--structures\index{$SU(m)$} by taking $\Psi_L=\Psi=\Psi_R$. It is instructive to work out the corresponding spinors $[\mc{A}(\Psi_L)\otimes\Psi_R]$ and $[\Psi_L\otimes\Psi_R]$ explicitly. For this, recall that an $SU(m)$--structure\index{$SU(m)$} on $W^{2m}$ can be equivalently defined in terms of a non--degenerate $2$--form $\omega$ and a decomposable complex $m$--form $\Omega$ satisfying certain algebraic relations (cf. for instance~\cite{hi97}) and the example of $SU(3)$\index{$SU(3)$} at the end of Section~\ref{spinstruc}). We then find for the corresponding straight structure the $W\oplus W^*$--spinors
$$
[\mc{A}(\Psi_L)\otimes\Psi_R]=(-1)^{m(m+1)/2}e^{-i\omega},\quad[\Psi_L\otimes\Psi_R]=\Omega.
$$
In particular, this example shows that not any generalised $SU(m)$--structure\index{generalised $SU(m)$--structure} arises as the $B$--field transform of a classical $SU(m)$--structure\index{$SU(m)$}.

\paragraph{Generalised $G_2$\index{generalised $G_2$--structure}--structures~\cite{wi06}.}
A {\em generalised $G_2$\index{$G_2$}--structure} is the choice of an embedding $G_2\times G_2\hookrightarrow Spin(7,0)\times Spin(0,7)$\index{$G_2$}.
We have $\cliff(7)\cong\End_{\R}(P_7)$ according to Table~(\ref{clifftable}), where $P_7\cong\R^8$ carries an invariant, positive definite inner product whose unit sphere is $S^7=Spin(7)/G_2$\index{$G_2$}\index{$Spin(7)$}. Any real unit spinor $\Psi$ thus induces a $G_2$\index{$G_2$}--structure. Since the complex spin representation is just $\Delta_7=P\otimes\C$, a pair of real unit spinors $(\Psi_L,\Psi_R)$ induces the $G_{2L}\times G_{2R}$--invariant spinors
$$
\rho_0=[\Psi_L\otimes\Psi_R]_B^{ev},\quad\rho_1=[\Psi_L\otimes\Psi_R]_B^{od}.
$$
Since
$$
\gtilde\bullet\rho_0=\rho_1,\quad\gtilde\bullet\rho_1=-\rho_0,
$$
a generalised $G_2$\index{generalised $G_2$--structure}--structure is characterised by either $[\Psi_L\otimes\Psi_R]_B^{ev}$ or $[\Psi_L\otimes\Psi_R]_B^{od}$. Again, this comprises the $B$--field transformation\index{$B$--field transformation} of a classical $G_2$\index{$G_2$}--structure. In the straight case, we find
$$
[\Psi\otimes\Psi]^{ev}=1-\star\varphi,\quad[\Psi\otimes\Psi]^{od}=-\varphi+vol_g,
$$
where $\varphi$ is the stable 3--form characterising the $G_2$\index{$G_2$}--structure on $W^7$.

\paragraph{Generalised $Spin(7)$\index{generalised $Spin(7)$--structure}--structures~\cite{wi06}.}
A {\em generalised $Spin(7)$\index{$Spin(7)$}--structure} is the choice of an embedding $Spin(7)\times Spin(7)\hookrightarrow Spin(8,0)\times Spin(0,8)$\index{$Spin(7)$}. Using again Table~(\ref{clifftable}), $\cliff(8)\cong\End(P_8)$, where $P_8$ splits into the two 8--dimensional representation spaces $S_{8\pm}$ which also carry an invariant, positive definite inner product and whose complexification gives the complex spin representation $\Delta_8$ of $Spin(8)$. The stabiliser of a real chiral unit spinor is isomorphic with $Spin(7)$\index{$Spin(7)$}. However, there are two conjugacy classes of $Spin(7)$\index{$Spin(7)$}, each of which stabilises a unit spinor in $S_{8+}$ or $S_{8-}$ respectively. Here, we will only consider the case where both spinors $\Psi_L$ and $\Psi_R$ are of equal chirality (this corresponds to the {\em even} type of~\cite{wi06}, while the case of opposite chirality leads to {\em odd} type structures). The $W\oplus W^*$--spinor $$
\rho=[\Psi_L\otimes\Psi_R]_B
$$ 
is even and invariant under $Spin(7)_L\times Spin(7)_R$\index{$Spin(7)$}. Moreover, $\rho$ is self--dual for $\gtilde$, i.e. 
$$
\gtilde\bullet\rho=\rho.
$$ 
In the straight case we obtain 
$$
[\Psi\otimes\Psi]=1-\Omega+vol_g,
$$
where $\Omega$ is the self--dual $4$--form characterising the classical $Spin(7)$\index{$Spin(7)$}--structure on $W^8$.
%
%
%
%
%
\section{The generalised tangent bundle\index{generalised tangent bundle}}
\label{gentangent}
In this section, we want to realise the first step of our $G$--structure approach to type II\index{type II} supergravity, namely to associate a $G$--structure with the field content $(g,H,\Psi_L,\Psi_R)$, where $g$ is a metric, $H$ a closed (not necessarily integral) $3$--form, and $\Psi_L$, $\Psi_R$ are two $TM$--spinors of unit length. 
%
%
\subsection{Twisting with an $H$--flux}
\label{gerbes}
The generalised setup discussed in the previous sections makes also sense globally. The bundle $TM\oplus T^*M$ carries a natural orientation and inner product provided by contraction, which extends the $SO(n,n)$--structure of every fibre $T_pM\oplus T^*_pM$. Further, $TM$ and $T^*M$ are isotropic subbundles. We can then speak about generalised Riemannian metrics, $SU(m)$--structures\index{$SU(m)$} etc. for $TM\oplus T^*M$. In addition, we also want to incorporate the $H$--flux, and for this, we take up the local approach emphasised in Section~\ref{gstructures}. There, we viewed a vector bundle as being defined by a collection of $G$--valued transition functions.

Assume $M^n$ to be oriented; the tangent bundle is then associated with some family of transition functions $s_{ab}\in GL(n)_+$. We can extend these to transition functions
$$
S_{ab}:U_a\cap U_b\to SO(n,n),\quad S_{ab}(p)=\left(\begin{array}{cc} s_{ab} & 0\\ 0 & s_{ab}^{-1tr}\end{array}\right)
$$
for $TM\oplus T^*M$. By refining the cover $\{U_a\}$ if necessary, we can assume it to be convex. The closed $3$--form $H$ is locally exact, that is, $H_{|U_a}=dB^{(a)}$ and we define
$$
\beta^{(ab)}=(B^{(a)}-B^{(b)})_{|U_a\cap U_b}\in\Omega^2(U_a\cap U_b).
$$
By design, the $2$--forms $\beta^{(ab)}$ are closed. Trivialising $TM$ over $U_a$ or $U_b$, we can think of these $2$--forms as maps $U_a\cap U_b\to\Lambda^2\R^n$. We indicate the trivialisation we use by the subscript $a$ or $b$, that is for $p\in U_a\cap U_b$,
$$
\beta^{(ab)}(p)=[a,p,\beta^{(ab)}_a(p)]=[b,p,\beta^{(ab)}_b(p)],
$$
so that $\beta^{(ab)}_a=s_{ab}^*\beta^{(ab)}_b$.
This gives rise to the twisted transition functions
$$
\sigma_{ab}=S_{ab}\circ e^{2\beta^{(ab)}_b}=e^{2\beta^{(ab)}_b}\circ S_{ab}.
$$
Indeed, whenever $U_a\cap U_b\cap U_c\not=0$, we have
\begin{eqnarray*}
\sigma_{ab}\circ \sigma_{bc} & = & S_{ab}\circ e^{2\beta^{(ab)}_b}\circ e^{2\beta^{(bc)}_b}\circ S_{bc}\\
& = & S_{ab}\circ e^{2(B^{(a)}_b-B^{(b)}_b+B^{(b)}_b-B^{(c)}_b)}\circ S_{bc}\\
& = & S_{ab}\circ S_{bc}\circ e^{2(B^{(a)}_c-B^{(c)}_c)}=S_{ac}\circ e^{2\beta^{(ac)}_c}\\
& = & \sigma_{ca}^{-1}.
\end{eqnarray*}
We define the {\em generalised tangent bundle}\index{generalised tangent bundle} by
\begin{equation}\nonumber
\mb{E}=\mb{E}(H)=\coprod U_a\times(\R^n\oplus\R^{n*})/\sim_{\sigma_{ab}}.
\end{equation}
Up to isomorphism, this bundle only depends on the closed $3$--form $H$. Indeed, assume we are given a different convex cover $\{U_a'\}$ together with locally defined $2$--forms $B^{(a)'}\in\Omega^2(U'_a)$ such that $H_{|U'_a}=dB^{(a)'}$, resulting in a new family of transition functions $\sigma'_{ab}=S_{ab}\circ\exp(2\beta'_b)$. Now on the intersection $V_a=U_a\cap U'_a$ we have $d(B^{(a)}_{|V_a}-B^{(a)'}_{|V_a})=0$, hence
\begin{equation}\nonumber
B^{(a)}_{|V_a}-B^{(a)'}_{|V_a}=dG^{(a)}\quad\mbox{for }G^{(a)}\in\Omega^1(V_a).
\end{equation}
One readily verifies the family $G_a=\exp(dG^{(a)}_a)$ to define a gauge transformation, i.e.
\begin{equation}\nonumber
\sigma'_{ab}=G^{-1}_a\circ \sigma_{ab}\circ G_b\mbox{ on }V_{ab}\not=\emptyset.
\end{equation}
In particular, the bundles defined by the families $\sigma_{ab}$ and $\sigma'_{ab}$ respectively are isomorphic, where the isomorphism is provided by
$$
[a,p,v]_{\sim_{s_{ab}}}\mapsto [a,p,G_a(v)]_{\sim_{s'_{ab}}}.
$$
Since the transition functions $\sigma_{ab}$ take values in $SO(n,n)_+$, the invariant orientation and inner product $(\cdot\,,\cdot)$ on $\R^n\oplus\R^{n*}$ make sense globally and turn $\mb{E}$ into an oriented pseudo--Riemannian vector bundle. 

Again, we can consider reductions inside this $SO(n,n)_+$--structure, for instance to $SO(n,0)\times SO(0,n)$.

\begin{definition}
A {\em generalised Riemannian metric} for the generalised tangent bundle\index{generalised tangent bundle} $\mb{E}(H)$ is a reduction from its structure group $SO(n,n)_+$ to $SO(n,0)\times SO(0,n)$.
\end{definition}

Since the group $SO(n,0)\times SO(0,n)$ preserves a decomposition of $\R^n\oplus\R^{n*}=V^+\oplus V^-$ into a positive and negative definite subspace $V^+$ and $V^-$, we can equivalently define a generalised Riemannian metric by the choice of a maximally positive definite subbundle $\mb{V}^+\subset\mb{E}(H)$. This bundle provides a splitting of the exact sequence

\begin{equation}\label{sequence}
0\to T^*M\stackrel{i}{\to} \mb{E}\stackrel{\pi}{\to} TM\to 0.
\end{equation}

Here, $i:T^*M\to \mb{E}$ is the canonical inclusion: Any differential form transforms under $s_{ab}^*$ and is acted on trivially by $\exp(\beta^{(ab)})$, therefore defining a section of $\mb{E}$. Identifying $T^*M$ with its image $i(T^*M)$ in $\mb{E}$, $T^*M\cap\mb{V}^+=\{0\}$, for $T^*M$ is isotropic. Hence the projection $\pi:\mb{E}\to TM$ restricted to $\mb{V}^+$ is injective, so that a generalised Riemannian structure defines a splitting of the exact sequence~(\ref{sequence}). We obtain a lift from $TM$ to $\mb{V}^+$ which in accordance with the notation used in the previous section we denote by $X^+$, for $X$ a vector field $X\in \Gamma(TM)$. Locally, $X$ corresponds to smooth maps $X_a:U_a\to\R^n$ such that $X_a=s_{ab}X_b$, while for $X^+_a:U_a\to\R^n\oplus\R^{n*}$, the relation $X^+_a=\sigma_{ab}X^+_b$ holds. As before, $\mb{V}^+$ is obtained as the graph of a linear isomorphism $P^+_a:\R^n\to\R^{n*}$, so
\begin{equation}\nonumber
X^+_a=X_a\oplus P^+_aX_a.
\end{equation}
From the transformation rule on $\{X^+_a\}$, we deduce 
\begin{equation}\label{genmetloc}
\beta^{(ab)}_a=P_a^+-s_{ab}^*P^+_bs_{ba}.
\end{equation}
The symmetric part $g_a=(P_a^++P_a^{+{\rm tr}})/2$ is therefore positive definite, and since $\beta^{(ab)}$ is skew--symmetric, the symmetrisation of the right hand side vanishes. Hence
\begin{equation}\nonumber
\frac{1}{2}(P_a^+-s_{ab}^*P^+_bs_{ba}+P_a^{+{\rm tr}}-s_{ab}^*P^{+{\rm tr}}_bs_{ba})=g_a-s_{ab}^*g_as_{ab}=0,
\end{equation}
so that the collection $g_a:U_a\to\odot^2\R^{n*}$ of positive definite symmetric $2$--tensors patches together to a globally defined metric. Conversely, a Riemannian metric $g$ induces a generalised Riemannian structure on $\mb{E}(H)$: The maps $P_a=B^{(a)}_a+g_a$ induce local lifts of $TM$ to $\mb{E}$ which give rise to a global splitting of~(\ref{sequence}).

\begin{proposition}
A generalised Riemannian structure is characterised by the datum $(g,H)$, where $g$ is a Riemannian metric and $H$ a closed $3$--form. 
\end{proposition}

\begin{remark}\hfill\newline
\noindent(i) From~(\ref{genmetloc}) we also conclude that the skew--symmetric part of $P_a$ is $B^{(a)}_a$. Hence we have local isomorphisms $V^+_{|U_a}=e^{2B^{(a)}}D^+_{|U_a}$, where $TM\oplus T^*M=D^+\oplus D^-$. In this way, we can think of the local model of a generalised Riemannian metric as a $B$--field transformed Riemannian metric. 

\noindent(ii) Of course, the negative definite subbundle $V^-$ also defines a splitting of~(\ref{sequence}). The lift of a vector field $X$ is then induced by $X^-_a=X_a\oplus P^-_aX_a$ with $P^-_a=-g_a+B^{(a)}_a$.
\end{remark}

For generalised $SU(m)$\index{$SU(m)$}--, $G_2$\index{$G_2$}-- or $Spin(7)$\index{$Spin(7)$}--structures, we need to speak about $\mb{E}$--spinor fields. This shall occupy us next.
%
%
\subsection{Spinors}
As discussed in Section~\ref{spinstruc}, we first need to exhibit a spin structure\index{spin structure}, i.e. a $Spin(n,n)_+$--valued family of functions $\widetilde{\sigma}_{ab}$ satisfying~(\ref{cocycle}) and covering $\sigma_{ab}$, that is $\pi_0\circ\widetilde{\sigma}_{ab}=\sigma_{ab}$. In the situation present, we can make a canonic choice: Exponentiating $\beta^{(ab)}$ to $Spin(n,n)_+$ and considering $GL(n)_+$ as its subgroup, we define
\begin{equation}\nonumber
\widetilde{\sigma}_{ab}=\widetilde{S}_{ab}\bullet e^{\beta^{(ab)}_b}=e^{\beta^{(ab)}_a}\bullet\widetilde{S}_{ab}.
\end{equation} 
The even and odd spinor bundles associated with $\mb{E}=\mb{E}(H)$ are
\begin{equation}\nonumber
\mb{S}(\mb{E})_{\pm}=\coprod\limits_aU_a\times S_{\pm}/\sim_{\widetilde{\sigma}_{ab}}.
\end{equation}
An $\mb{E}$--spinor field $\rho$ is thus represented by a collection of smooth maps $\rho_a:U_a\to S_{\pm}$ with $\rho_a=\widetilde{\sigma}_{ab}\bullet\rho_b$. 
Since as a {\em vector space}, $S_{\pm}=\Lambda^*\R^{*n}$, it is tempting to think of chiral $\mb{E}$--spinors as even or odd differential forms. As we have already remarked (cf.~(\ref{spinGL})), $s_{ab}\in GL(n)$ acts on $\rho_b$ via $\sqrt{\det s_{ab}}\bullet s_{ab}^*\rho_b$, where $s_{ab}^*$ denotes the induced action of $GL(n)$ on forms, so that the family $\{\rho_a\}$ does not transform as a differential form. To remedy this, we pick a nowhere vanishing $n$--vector field $\nu$, that is, a collection of smooth maps $\nu_a:U_a\to\Lambda^n\R^{n*}$, $\nu_a=\lambda^{-2}_a\nu_0$, where the coefficient of $\nu_a$ in $C^{\infty}(U_a)$ is assumed to be strictly positive. The notation $\lambda^{-2}_a$ is introduced to ease notation in the subsequent computations. Since $\nu$ is globally defined, the coefficients transform under $\lambda_a^{-2}=\det s_{ab}\cdot \lambda_b^{-2}$. We then define an isomorphism $\mf{L}^{\nu}:\Gamma\big(\mb{S}(E)_{\pm}\big)\to\Omega^{ev,od}(M)$ by
$$
\mf{L}^{\nu}_a:(\rho_a:U_a\to S_{\pm})\mapsto (e^{-B^{(a)}_a}\wedge\lambda_a\cdot\rho_a:U_a\to\Lambda^{ev,od}\R^{n*}).
$$
We need to show that this transforms correctly under the action of the transition functions $s_{ab}$ on $\Lambda^{ev,od}T^*M$, using the fact that $\rho_a=\widetilde{\sigma}_{ab}\bullet\rho_b$. Indeed, we have on $U_a\cap U_b\not=\emptyset$
\begin{eqnarray*}
s_{ab}^*(e^{-B^{(b)}_b}\wedge\lambda_b\cdot\rho_b) & = & \lambda_a\sqrt{\det s_{ab}}\cdot e^{-B^{(b)}_a}\wedge s^*_{ab}\rho_b\nonumber\\
& = & \lambda_a\cdot e^{-B^{(a)}_a}\wedge e^{\beta^{(ab)}_a}\wedge \sqrt{\det s_{ab}}\cdot s^*_{ab}\rho_b\nonumber\\
& = & \lambda_a\cdot e^{-B^{(a)}_a}\wedge\widetilde{\sigma}_{ab}\bullet\rho_b\\
& = & \lambda_a\cdot e^{-B^{(a)}_a}\wedge\rho_a,
\end{eqnarray*}
or equivalently, $\mf{L}^{\nu}_a\circ\widetilde{\sigma}_{ab}=s_{ab}^*\circ \mf{L}^{\nu}_b$. In the same vein, the $Spin(n,n)$--invariant form $\langle\cdot\,,\cdot\rangle$ induces as above a globally defined inner product on $\Gamma(\mb{S})$ by
\begin{equation}\nonumber
\langle\rho,\tau\rangle=\nu\big([\mf{L}^{\nu}(\rho)\wedge\widehat{\mf{L}^{\nu}(\tau)}]^n\big),
\end{equation}
where $\nu$ is a nowhere vanishing $n$--vector field. 

The fact that the generalised tangent bundle\index{generalised tangent bundle} is obtained out of twisting the transition functions of $TM\oplus T^*M$ with {\em closed} $B$--fields bears an important consequence. Given the choice of $\nu$ with local coefficients $\lambda_a^{-2}$, we define a parity reversing map \begin{equation}\label{spinder}
d_{\nu}:\Gamma\big(\mb{S}(E)_{\pm}\big)\to\Gamma\big(\mb{S}(E)_{\mp}\big),\quad(d_{\nu}\rho)_a=\lambda^{-1}_a\cdot d_a(\lambda_a\cdot\rho_a),
\end{equation}
where $d_a$ is the usual differential applied to forms $U_a\to\Lambda^*\R^{n*}$, i.e. $(d\alpha)_a=d_a\alpha_a$. This definition gives indeed rise to an $\mb{E}$--spinor field, for
\begin{eqnarray*}
\widetilde{\sigma}_{ab}\bullet(d_{\nu}\rho)_b & = & \lambda^{-1}_b\cdot e^{\beta^{(ab)}_a}\wedge \sqrt{\det s_{ab}}\cdot s_{ab}^*d_b(\lambda_b\cdot\rho_b)\\
& = & \lambda_a^{-1}\cdot d_a\big(e^{\beta^{(ab)}_a}\wedge s_{ab}^*(\lambda_b\cdot\rho_b)\big)\\
& = & \lambda_a^{-1}\cdot d_a(\lambda_a\cdot e^{\beta^{(ab)}_a}\wedge\sqrt{\det s_{ab}}\cdot s^*_{ab}\rho_b)\\
& = & \lambda^{-1}_a\cdot d_a(\lambda_a\cdot\rho_a)\\
& = & (d_{\nu}\rho)_a.
\end{eqnarray*}
The operator $d_{\nu}$ squares to zero and therefore induces an elliptic complex on $\Gamma\big(\mb{S}(E)_{\pm}\big)$. As a corollary of the next proposition, we deduce that this elliptic complex actually computes the so--called {\em twisted cohomology}, where on replaces the usual differential $d$ of de Rham cohomology by the twisted differential $d_H=d+H\wedge$.

\begin{proposition}\label{twisteddiff}
Let $\rho\in\Gamma\big(\mb{S}(E)\big)$. Then
\begin{equation}\nonumber
\mf{L}^{\nu}(d_{\nu}\rho)=d_H\mf{L}^{\nu}(\rho).
\end{equation}
\end{proposition}

\begin{proof}
This follows from a straightforward local computation:
\begin{eqnarray*}
\big(\mf{L}^{\nu}(d_{\nu}\rho)_a\big) & = & \mf{L}^{\nu}_a(d_{\nu}\rho)_a\\
& = & e^{-B^{(a)}_a}\wedge\lambda_a\cdot(d_{\nu}\rho)_a\\
& = & e^{-B^{(a)}_a}\wedge d_a(\lambda_a\cdot\rho_a)\\
& = & d_{H_a}(e^{-B^{(a)}_a}\wedge\lambda_a\cdot\rho_a)\\
& = & d_{H_a}\mf{L}^{\nu}_a(\rho_a).
\end{eqnarray*}
\end{proof}

In presence of a generalised Riemannian metric, we can make a canonic choice for $\nu$, namely we pick the dual $\nu_g$ of the Riemannian volume form $vol_g$. In this case, we write $\mf{L}^{\nu_g}=\mf{L}$ and $d_{\nu_g}=d$. Moreover, we obtain again an operator $\gtilde=vol_{\mb{V}^-}\bullet$ for which we find as above:

\begin{proposition}
The action of $\widetilde{\mc{G}}=vol_{V^-}\bullet$ on $\mb{S}_{\pm}(\mb{E})$ is given by
\begin{equation}\nonumber
\mf{L}(\gtilde\bullet \rho)=\left\{\begin{array}{ll}$n$\mbox{ even:}\quad&\star_g\widehat{\mf{L}(\rho)}\\$n$\mbox{ odd:}\quad&\star_g\widehat{\widetilde{\mf{L}(\rho)}}\end{array}\right.,
\end{equation}
where $g$ is the Riemannian metric associated with $\mb{V}^+$.
\end{proposition}

If the underlying manifold is spinnable, we can again identify bispinors with $\mb{E}$--spinors via the map
$$
[\cdot,\cdot]^{\mc{G}\,ev,od}:\Gamma\big(\Delta_n(TM)\otimes\Delta_n(TM)\big)\stackrel{[\cdot,\cdot\,]}{\longrightarrow}\Omega^{ev,od}(M)\otimes\C\stackrel{\mf{L}^{-1}}{\longrightarrow}\Gamma\big(\mb{S}(\mb{E})_{\pm}\otimes\C\big).
$$
A vector field $X\in\Gamma(TM)$ acts on $TM$--spinor fields via the inclusion $TM\hookrightarrow\cliff(TM,g)$ and Clifford multiplication\index{Clifford multiplication}. On the other hand, we can lift $X$ to sections $X^{\pm}$ of $\mb{V}^{\pm}$ which act on $\mb{E}$--spinor fields via the inclusion $\mb{V}^{\pm}\hookrightarrow\cliff(\mb{E})$ and Clifford multiplication\index{Clifford multiplication}. As in the previous section, we find that these actions are compatible in the following sense.

\begin{proposition}
We have
\begin{eqnarray*}\nonumber
\,[X\cdot\Psi_L\otimes\Psi_R]^{\mc{G}} & = & (-1)^{n(n-1)/2}X^+\bullet[\Psi_L\otimes\Psi_R]^{\mc{G}},\\
\,[\Psi_L\otimes Y\cdot\Psi_R]^{\mc{G}} & = & -Y^-\bullet\widetilde{[\Psi_L\otimes\Psi_R]^{\mc{G}}}.
\end{eqnarray*}
\end{proposition}

In particular, we get as an

\begin{corollary}\label{selfdual}
Let $\Psi_{L,R}\in\Delta_n$. 

\noindent{\rm (i)} If $n=2m$ and $\Psi_{L,R}$ are chiral, then for $\Psi_R\in\Delta_{\pm}$
\begin{equation}\nonumber
\gtilde\bullet[\Psi_L\otimes\Psi_R]^{\mc{G}}=\pm (-1)^{m(m-1)/2}i^m[\Psi_L\otimes\Psi_R]^{\mc{G}}.
\end{equation}
\noindent{\rm (ii)} If $n=2m+1$, then
\begin{equation}\nonumber
\gtilde\bullet[\Psi_L\otimes\Psi_R]^{\mc{G}}=(-1)^{m(m-1)/2}i^{m+1}\widetilde{[\Psi_L\otimes\Psi_R]^{\mc{G}}}.
\end{equation}
\end{corollary}

Now it is clear how we can describe reductions to $SU(m)\times SU(m)$\index{$SU(m)$}, $G_2\times G_2$\index{$G_2$} or $Spin(7)\times Spin(7)$\index{$Spin(7)$}.

\begin{example}
We consider a {\em generalised $SU(3)$\index{$SU(3)$}--structure}.
In dimension $6$, a reduction from $Spin(6)$ to $SU(3)$\index{$SU(3)$} is induced by a unit spinor, and there is no obstruction against existence as we saw in Section~\ref{spinstruc}. A reduction from the generalised Riemannian metric structure given by $(g,H)$ to $SU(3)\times SU(3)$\index{$SU(3)$} can be therefore characterised in terms of two unit spinors $(\Psi_L,\Psi_R)$, giving rise to the $\mb{E}$--spinors
$$
\rho_0=[\mc{A}(\Psi_L)\otimes\Psi_R]^{\mc{G}},\quad\rho_1=[\Psi_L\otimes\Psi_R]^{\mc{G}}.
$$
The corresponding {\em differential forms} are just
$$
\mf{L}(\rho_0)=e^{-i\omega},\quad\mf{L}(\rho_1)=\Omega
$$
(cf. Section~\ref{specorb}).
\end{example}
%
%
%
%
%
\section{The field equations}
\label{fieldequs}
The second step of the $G$--structure ansatz consists in recovering the type II\index{type II} field equations~(\ref{mastereqgrav})and~(\ref{mastereqdil}) by an integrability condition on the algebraic objects defining the $G$--structure. Here, we will deal with the case $F=0$. For a treatment with non--trivial R-R--fields\index{Ramond--Ramond field} in the case of $SU(3)$\index{$SU(3)$}-- and $G_2$\index{$G_2$}--structures, see~\cite{jewi05b}. Throughout this section, $\nabla$ denotes the Levi--Civita connection. By a generalised $G$--structure, we shall mean a generalised structure characterised by a collection of decomposable bispinors such as generalised $SU(m)$\index{$SU(m)$}--, $G_2$\index{$G_2$}-- or $Spin(7)$\index{$Spin(7)$}--structures.
%
%
\subsection{Integrable\index{integrable} generalised $G$--structures}
\label{intgenstruc} 
\begin{definition}\label{formpicture}
Let $\{\rho_0,\ldots,\rho_l\}$ be a collection of $\mb{E}$--spinors defining a generalised $G$--structure. Then this structure is called {\em integrable}\index{integrable} if 
\begin{equation}\label{intcond}
d_{\nu}\rho_i=0,\quad d_{\nu}\gtilde\bullet\rho_i=0,\quad i=0,\ldots,l
\end{equation}
for some $\nu\in\Gamma(\Lambda^nTM)$
\end{definition}

\begin{example}
Consider a Riemannian manifold $(M^{2m},g)$ whose holonomy is contained in $SU(m)$\index{$SU(m)$}. As a consequence, $M$ carries an $SU(m)$\index{$SU(m)$}--structure with $d\omega=0$ and $d\Omega=0$. Then the corresponding straight structure characterised by $\rho_0$ and $\rho_1$ is also integrable\index{integrable}, for $d\rho_{0,1}=0$ if and only if 
$$
d[\mc{A}(\Psi)\otimes\Psi]=(-1)^{m(m+1)/2}de^{-i\omega}=0,\quad d[\Psi\otimes\Psi]=d\Omega=0
$$
by Proposition~\ref{twisteddiff}. From this point of view, the integrability condition~(\ref{intcond}) generalises the holonomy condition of classical $SU(m)$\index{$SU(m)$}--structures.
\end{example}

If the generalised $G$--structure induces a Riemannian metric, then we can define the {\em dilaton field}\index{dilaton field} $\phi\in C^{\infty}(M)$ via
$$
\nu=e^{2\phi}\cdot\nu_g.
$$
We can then write~(\ref{intcond}) as a form equation
$$
d_He^{-\phi}\mf{L}(\rho_i)=0,\quad \pm d_He^{-\phi}\star\widehat{\mf{L}(\rho_i)}=0,\quad i=0,\ldots,l.
$$

\begin{remark}
The appearance of the dilaton field\index{dilaton field} may seem artificial. However, there are two reasons to it: Firstly, we will prove a no--go theorem in the next section which asserts that a constant dilaton field\index{dilaton field} implies $H=0$ if the structure is integrable\index{integrable}. In conjunction with the theorem we are going to prove in a moment, this means that the only integrable\index{integrable} generalised structures which occur are straight structures. Secondly, we can derive the integrability condition~(\ref{intcond}) on generalised $SU(3)$\index{$SU(3)$}-- and $G_2$\index{$G_2$}--structures from Hitchin's variational principle which requires an a priori identification of $\mb{E}$--spinors with forms and thus the choice of some dilaton field\index{dilaton field}~\cite{jewi05b},~\cite{wi06}.
\end{remark}

The following theorem links integrability of $SU(m)$\index{$SU(m)$}--structures into the supersymmetry equations~(\ref{susyeq1}) and~(\ref{susyeq2}). For generalised $G_2$\index{generalised $G_2$--structure}-- and $Spin(7)$\index{generalised $Spin(7)$--structure}--structures, see~\cite{wi06}.

\begin{theorem}\label{integrability}
Let $(\rho_0,\rho_1)=([\mc{A}(\Psi_L)\otimes\Psi_R]^{\mc{G}},[\Psi_L\otimes\Psi_R]^{\mc{G}})$ be a generalised $SU(m)$\index{$SU(m)$}--structure and $\phi\in C^{\infty}(M)$. Then
$$
d_He^{-\phi}[\mc{A}(\Psi_L)\otimes\Psi_R]=0,\quad d_He^{-\phi}[\Psi_L\otimes\Psi_R]=0
$$
holds, i.e. the generalised $SU(m)$\index{$SU(m)$}--structure is integrable\index{integrable}, if and only if the equations
\begin{equation}\label{spinorfield}
\begin{array}{ll}
\nabla_X\Psi_L-\frac{1}{4}(X\llcorner H)\cdot\Psi_L=0,\quad& (d\phi-\frac{1}{2}H)\cdot\Psi_L=0,\\
\nabla_X\Psi_R+\frac{1}{4}(X\llcorner H)\cdot\Psi_R=0,\quad& (d\phi+\frac{1}{2}H)\cdot\Psi_R=0,
\end{array}
\end{equation}
hold.
\end{theorem}

\begin{proof}
For computational purposes, it will be convenient to consider the equation 
$$
d[\Psi_L\otimes\Psi_R]=-(\alpha+H)\wedge[\Psi_L\otimes\Psi_R],\; d[\mc{A}(\Psi_L)\otimes\Psi_R]=-(\alpha+H)\wedge [\mc{A}(\Psi_L)\otimes\Psi_R]
$$
for a $1$--form $\alpha$ instead of the dilaton. The key for solving this set of equations is the decomposability of the spinor: it makes the spinor and its associated differential form ``self--dual'' in the sense of Corollary~\ref{selfdual}. From~(\ref{selfdualgensum}) $\star\widehat{[\Psi_L\otimes\Psi_R]}=(-1)^{m(m-1)/2}i^m[\Psi_L\otimes\Psi_R]$, so that 
\begin{eqnarray}
d\star\widehat{[\Psi_L\otimes\Psi_R]} & = & d\big((-1)^{m(m-1)/2}i^m[\Psi_L\otimes\Psi_R]\big)\nonumber\\
& = & -(-1)^{m(m-1)/2}i^m(\alpha+H)\wedge[\Psi_L\otimes\Psi_R]
\end{eqnarray} 
and similarly for $[\mc{A}(\Psi_L)\otimes\Psi_R]$. We recall that 
\begin{equation}\nonumber
\widetilde{[\mc{A}(\Psi_L)\otimes\Psi_R]}=[\mc{A}(\Psi_L)\otimes\Psi_R],\quad\widetilde{[\Psi_L\otimes\Psi_R]}=(-1)^m[\Psi_L\otimes\Psi_R],
\end{equation} 
as well as the general rules for forms $R\in\Omega^*(M^{2m})$, namely 
\begin{equation}\nonumber
\star\widehat{R}=(-1)^m\widehat{\widetilde{\star R}},\quad\star(\alpha\wedge R)=\alpha\llcorner\star\widetilde{R},\quad\widehat{dR}=-d\widehat{\widetilde{R}},\quad d^*R=-\star d\star R.
\end{equation} 
From these we deduce
\begin{eqnarray*}
d^*[\Psi_L\otimes\Psi_R] & = & -(\alpha+H)\llcorner\widetilde{[\Psi_L\otimes\Psi_R]},\\ 
d^*[\mc{A}(\Psi_L)\otimes\Psi_R] & = & \phantom{-}(\alpha+H)\llcorner [\mc{A}(\Psi_L)\otimes\Psi_R].
\end{eqnarray*}
From Proposition~\ref{commutation} follows immediately a technical lemma we need next.

\begin{lemma}
Let $\alpha$ be a $1$--form. Its metric dual will be also denoted by $\alpha$. Then
\begin{eqnarray*}
\alpha\wedge[\Psi_1\otimes\Psi_2] & = & \phantom{-}\frac{1}{2}\big((-1)^m[\alpha\cdot\Psi_1\otimes\Psi_2]-\widetilde{[\Psi_1\otimes \alpha\cdot\Psi_2]}\big)\\
\alpha\llcorner[\Psi_1\otimes\Psi_2] & = & -\frac{1}{2}\big((-1)^m[\alpha\cdot\Psi_1\otimes\Psi_2]+\widetilde{[\Psi_1\otimes \alpha\cdot\Psi_2]}\big)
\end{eqnarray*}
In particular, if $e_k$ defines a local orthonormal basis, then for a $3$--form $H$ we find
\begin{eqnarray*}
H\wedge[\Psi_1\otimes\Psi_2] & = & \phantom{-}\frac{(-1)^m}{8}\big[H\cdot\Psi_1\otimes\Psi_2-\sum_ke_k\cdot\Psi_1\otimes(e_k\llcorner H)\cdot\Psi_2\big]\\
& & +\frac{1}{8}\big[\Psi_1\otimes H\cdot\Psi_2-\sum_k(e_k\llcorner H)\cdot\Psi_1\otimes e_k\cdot\Psi_2\big]^{\sim}\\
H\llcorner[\Psi_1\otimes\Psi_2] & = &  \phantom{-}\frac{(-1)^m}{8}\big[-H\cdot\Psi_1\otimes\Psi_2+\sum_ke_k\cdot\Psi_1\otimes(e_k\llcorner H)\cdot\Psi_2\big]\\
& & +\frac{1}{8}\big[\Psi_1\otimes H\cdot\Psi_2-\sum_k(e_k\llcorner H)\cdot\Psi_1\otimes e_k\cdot\Psi_2\big]^{\sim}.
\end{eqnarray*}
Moreover, since $d=\sum e_k\wedge\nabla_{e_k}$ and $d^*=-\sum e_k\llcorner\nabla_{e_k}$, we get
\begin{eqnarray*}
d[\Psi_1\otimes\Psi_2] & = & \frac{1}{2}\big((-1)^m\widetilde{[\mc{D}(\Psi_1\otimes\Psi_2)]}-\widetilde{[\widetilde{\mc{D}}(\Psi_1\otimes\Psi_2)]}\big),\\
d^*[\Psi_1\otimes\Psi_2] & = & \frac{1}{2}\big((-1)^m\widetilde{[\mc{D}(\Psi_1\otimes\Psi_2)]}+\widetilde{[\widetilde{\mc{D}}(\Psi_1\otimes\Psi_2)]}\big),
\end{eqnarray*}
with the {\em twisted Dirac operators\index{Dirac operator}} $\mc{D}$ and $\widetilde{\mc{D}}$ on $\Gamma\big(\Delta(T)\otimes\Delta(T)\big)$,
given locally by
\begin{eqnarray*}
\mc{D}(\Psi_1\otimes\Psi_2) & = & \sum
e_k\cdot\nabla_{e_k}\Psi_1\otimes\Psi_2+e_k\cdot\Psi_1\otimes\nabla_{e_k}\Psi_2\\
& = & \mb{D}\Psi_1\otimes\Psi_2+\sum e_k\cdot\Psi_1\otimes\nabla_{e_k}\Psi_2,\\
\widetilde{\mc{D}}(\Psi_1\otimes\Psi_2) & = & \sum \nabla_{e_k}\Psi_1\otimes
e_k\cdot\Psi_2+\Psi_1\otimes e_k\cdot\nabla_{e_k}\Psi_2\\ 
& = & \sum\nabla_{e_k}\Psi_1\otimes e_k\cdot\Psi_2+\Psi_1\otimes \mb{D}\Psi_2.
\end{eqnarray*}
\end{lemma}
Note that $\nabla$ and $[\cdot\,,\cdot]$ commute, since $\nabla$ is metric. As a result,
\begin{eqnarray}
\big[\mc{D}(\Psi_L\otimes\Psi_R)\big] & = & \phantom{-}(\alpha+H)\wedge[\Psi_L\otimes\Psi_R]+(-1)^m(\alpha+H)\llcorner[\Psi_L\otimes\Psi_R]\nonumber
\\
\big[\mc{D}\big(\mc{A}(\Psi_L)\otimes\Psi_R\big)\big] & = & \phantom{-}(-1)^m\big((\alpha+H)\wedge-(\alpha+H)\llcorner\big)[\mc{A}(\Psi_L)\otimes\Psi_R]\nonumber\\
\big[\widetilde{\mc{D}}(\Psi_L\otimes\Psi_R)\big] & = & -(-1)^m(\alpha+H)\wedge[\Psi_L\otimes\Psi_R]+(\alpha+H)\llcorner[\Psi_L\otimes\Psi_R]\nonumber\\
\big[\widetilde{\mc{D}}(\mc{A}(\Psi_L)\otimes\Psi_R)\big] & = & -\big((\alpha+H)\wedge+(\alpha+H)\llcorner\big)[\mc{A}(\Psi_L)\otimes\Psi_R].\label{twisteddiracev}
\end{eqnarray}
Using the previous lemma to compute the action of $\alpha+H$ on $[\Psi_L\otimes\Psi_R]$ the two first equations of~(\ref{twisteddiracev}) become
\begin{eqnarray}
\mb{D}\Psi_L\otimes\Psi_R+\sum e_k\cdot\Psi_L\otimes\nabla_{e_k}\Psi_R & = & -\alpha\cdot\Psi_L\otimes\Psi_R\nonumber\\
& & +\frac{1}{4}\sum e_k\cdot\Psi_L\otimes(e_k\llcorner H)\cdot\Psi_R\nonumber\\
& & -\frac{1}{4}H\cdot\Psi_L\otimes\Psi_R\label{tweq1}\\
\mb{D}\mc{A}(\Psi_L)\otimes\Psi_R+\sum e_k\cdot\mc{A}(\Psi_L)\otimes\nabla_{e_k}\Psi_R & = & -\alpha\cdot\mc{A}(\Psi_L)\otimes\Psi_R\nonumber\\
& &
+\frac{1}{4}\sum e_k\cdot\mc{A}(\Psi_L)\otimes(e_k\llcorner H)\cdot\Psi_R\nonumber\\
& &  -\frac{1}{4}H\cdot\mc{A}(\Psi_L)\otimes\Psi_R.\label{tweq2}
\end{eqnarray}
Contracting~(\ref{tweq1}) from the left hand side with $q(e_j\cdot\Psi_L,\,\cdot)$ gives
\begin{eqnarray*}
0 & = & \phantom{+}q\big(e_j\cdot\Psi_L,\mb{D}\Psi_L+\alpha\cdot\Psi_L+\frac{1}{4}H\cdot\Psi_L\big)\Psi_R\\
& & +\sum q(e_j\cdot\Psi_L,e_k\cdot\Psi_L)(\nabla_{e_k}-e_k\llcorner H)\cdot\Psi_R.
\end{eqnarray*}
We apply the conjugate linear operator $\mc{A}$ which commutes with the Levi--Civita connection $\nabla$ since it is $Spin(n)$--invariant, and get
\begin{eqnarray}
0 & = & \phantom{+}\overline{q\big(e_j\cdot\Psi_L,\mb{D}\Psi_L+\alpha\cdot\Psi_L+\frac{1}{4}H\cdot\Psi_L\big)}\mc{A}(\Psi_R)\nonumber\\
& & +\sum_k \overline{q(e_j\cdot\Psi_L,e_k\cdot\Psi_L)}(\nabla_{e_k}-\frac{1}{4}e_k\llcorner H\cdot)\mc{A}(\Psi_R)\label{coneq1}
\end{eqnarray}
Moreover, applying $\mc{A}\otimes \mc{A}$ to~(\ref{tweq2}) yields
\begin{eqnarray*}
\mb{D}\Psi_L\otimes\mc{A}(\Psi_R)+\sum e_k\cdot\Psi_L\otimes\nabla_{e_k}\Psi_R & = & \phantom{+}\frac{1}{4}\sum e_k\cdot\Psi_L\otimes (e_k\llcorner H)\cdot\mc{A}(\Psi_R)\\
& & -\frac{1}{4}H\cdot\Psi_l\otimes\mc{A}(\Psi_R)
\end{eqnarray*}
as
\begin{equation}\nonumber
[\mc{A}(\Psi)\otimes\mc{A}(\Phi)]=\widehat{m}\sum_Kq\big(\Psi,e_K\cdot\mc{A}(\Phi)\big)e_K=\widehat{m}\widehat{\overline{[\Phi\otimes\Psi]}}=\left\{\begin{array}{ll}m\mbox{ even:}\quad&\overline{[\Psi\otimes\Phi]}\\
m\mbox{ odd:}\quad&\widetilde{\overline{[\Psi\otimes\Phi]}}\end{array}\right.
\end{equation}
Contracting again with $q(e_j\cdot\Psi_L,\,\cdot)$ on the left gives 
\begin{eqnarray}
0 & = & \phantom{+}q(e_j\cdot\Psi_L,\mb{D}\Psi_L+\alpha\cdot\Psi_L+\frac{1}{4}H\cdot\Psi_L)\mc{A}(\Psi_R)\nonumber\\
& & +\sum_k q(e_j\cdot\Psi_L,e_k\cdot\Psi_L)\big(\nabla_{e_k}-\frac{1}{4}e_k\llcorner H\cdot\big)\mc{A}(\Psi_R).\label{coneq2}
\end{eqnarray}
Adding~(\ref{coneq1}) and~(\ref{coneq2}) yields
\begin{eqnarray}
0 & = & \phantom{+}{\rm Re}\, q(e_j\cdot\Psi_L,\mb{D}\Psi_L-\alpha\cdot\Psi_L-\frac{1}{4}H\cdot\Psi_L)\mc{A}(\Psi_R)\nonumber\\
& & +\sum {\rm Re}\, q(e_j\cdot\Psi_L,e_k\cdot\Psi_L)(\nabla_{e_k}-\frac{1}{4}e_k\llcorner H\cdot)\mc{A}(\Psi_R).
\end{eqnarray}
Now the real part of $q(e_j\cdot\Psi_L,e_k\cdot\Psi_L)$ vanishes unless $j=k$ when it equals $1$. This implies 
\begin{eqnarray*}
\nabla_{e_j}\mc{A}(\Psi_R) & = & \phantom{-}\frac{1}{4}(e_j\llcorner H)\cdot\mc{A}(\Psi_R)-\\
& & {\rm Re}\, q\big(e_j\cdot\Psi_L,\mb{D}\Psi_L+\alpha\cdot\Psi_L+\frac{1}{4}H\cdot\Psi_L\big)\mc{A}(\Psi_R),
\end{eqnarray*}
so that
\begin{eqnarray*}
\nabla_X\Psi_R & = & \frac{1}{4}(X\llcorner H)\cdot\Psi_R-{\rm Re}\, q\big(X\cdot\Psi_L,\mb{D}\Psi_L+\alpha\cdot\Psi_L+\frac{1}{4}H\cdot\Psi_L\big)\Psi_R
\end{eqnarray*}
We contract with $q(\mc{A}(\Psi_R),\cdot)$ to see that ${\rm Re}\,q(e_j\cdot\Psi_L,\ldots)=0$ as the remaining terms are purely imaginary, hence (upon applying $\mc{A}$)
\begin{equation}\nonumber
\nabla_X\Psi_R+\frac{1}{4}(X\llcorner H)\cdot\Psi_R=0.
\end{equation}
Using the second set of equations in~(\ref{twisteddiracev}) gives
\begin{equation}\nonumber
\nabla_X\Psi_L-\frac{1}{4}(X\llcorner H)\cdot\Psi_L=0.
\end{equation}
As a result, we find $\mb{D}\Psi_{R,L}=\mp3H/4\cdot\Psi_{R,L}$. On the other hand,~(\ref{tweq1}) now reads
\begin{equation}\nonumber
(\mb{D}\Psi_L-\alpha\cdot\Psi_L-\frac{1}{4}H\cdot\Psi_L)\otimes\Psi_R=0,
\end{equation}
so by contracting with $q(\Psi_R,\,\cdot)$ from the right hand side we obtain
\begin{equation}\nonumber
\mb{D}\Psi_L-\alpha\cdot\Psi_L-\frac{1}{4}H\cdot\Psi_L=0,
\end{equation}
which by the previous yields
\begin{equation}\nonumber
(\alpha-\frac{1}{2}H)\cdot\Psi_L=0.
\end{equation}
Similarly, we obtain $(\alpha+\frac{1}{2}H)\cdot\Psi_R=0$.
\end{proof}

\begin{remark}
If $H=0$, then $\phi\equiv const$, and we get two parallel spinors for the Levi--Civita connection, leaving us with two possibilities: Either the spinors coincide at one and thus at any point, or the two spinors are linearly independent everywhere, in which case we may assume that they are orthogonal. In either scenario, the holonomy reduces to the intersection of the stabilisers of $\Psi_L$ and $\Psi_R$ inside $Spin(n)$, giving rise to a well--defined classical $G$--structure. We therefore refer to these solutions as {\em classical}.
\end{remark}
%
%
\subsection{Geometric properties}
We now study some geometric properties of integrable\index{integrable} $G$--structures by using the formulation given in~(\ref{spinorfield}). In fact, most statements are valid for geometries defined by one parallel spinor, i.e. we suppose to be given a solution $\Psi$ to 
\begin{equation}\label{gravdil}
\nabla^H_X\Psi=\nabla_X\Psi+\frac{1}{4}(X\llcorner H)\cdot\Psi=0,\quad(d\phi+\frac{1}{2}H)\cdot\Psi=0.
\end{equation}
The key assumption here is that $H$ is closed, as we will presently see. 

To start with, we first compute the Ricci tensor. By results of~\cite{friv02}, the Ricci endomorphism $\ric^H$ of $\nabla^H$ with $H$ closed is given by
$$
\ric^H(X)\cdot\Psi=(\nabla^H_XH)\cdot\Psi,
$$
and relates to the metric Ricci tensor through
\begin{equation}\label{ricci}
\ric(X,Y)=\ric^H(X,Y)+\frac{1}{2}d^*H(X,Y)+\frac{1}{4}g(X\llcorner
H,Y\llcorner H).
\end{equation}
Consequently, the scalar curvature $S$ of the Levi--Civita connection is
$$
S=S^H+\frac{3}{2}\eunorm{H}^2,
$$
where $S^H$ is the scalar curvature associated with $\nabla^H$. Since $\Psi$ is parallel with respect to $\nabla^H$, the dilatino equation\index{dilatino equation} implies
$$
\ric^H(X)\cdot\Psi=\nabla^H_X(H\cdot\Psi)=-2(\nabla_X^Hd\phi)\cdot\Psi,
$$
hence $\ric^H(X)=-2\nabla_X^Hd\phi$. Now pick a frame that satisfies $\nabla_{e_i}e_j=0$ at a fixed
point, or equivalently,
$\nabla_{e_j}^He_k=e_k\llcorner e_j\llcorner H/2$. As the connection $\nabla^H$ is metric, we obtain
\begin{eqnarray*}
\ric^H(e_j,e_k) & = &
-2g(\nabla_{e_j}^Hd\phi,e_k)\\
& = & -2e_j.g(d\phi,e_k)+g(d\phi,\nabla^H_{e_j}e_k)\\
& = & -2e_j.e_k.\phi+e_k\llcorner e_j\llcorner H.\phi/2.
\end{eqnarray*}
The first summand is minus twice $\mc{H}^{\phi}$, the Hessian of $\phi$ evaluated in the basis $\{e_k\}$. Consequently, $\ric^H(X,Y)=-2\mc{H}^{\phi}(X,Y)-X\llcorner Y\llcorner H/2$, hence $S^H=2\Delta\phi$, where $\Delta(\cdot)=-{\rm Tr_g}\mc{H}^{(\cdot)}$ is the Riemannian Laplacian. In the situation where we have two spinors $\Psi_{L,R}$ parallel with respect to the connections $\nabla^{\pm H}$, we obtain from~(\ref{ricci})
$$
\ric(X,Y)=\frac{1}{2}\big(\ric^H(X,Y)+\ric^{-H}(X,Y)\big)+\frac{1}{4}g(X\llcorner H,Y\llcorner H),
$$ 
and thus the

\begin{theorem}\label{riccitensor}
The Ricci tensor $\ric$ and the scalar curvature $S$ of a metric of an integrable\index{integrable} generalised $G$--structure are given by
$$
\ric(X,Y)=-2\mc{H}^{\phi}(X,Y)+\frac{1}{2}g(X\llcorner H,
Y\llcorner H),\quad S=2\Delta
\phi+\frac{3}{2}\eunorm{H}^2.
$$
\end{theorem}

\begin{remark}
Note that Theorem~\ref{riccitensor} corrects an error in~\cite{wi06}, where in Theorem 4.9 and Proposition 5.7 the scalar curvature was stated to be $S=2\Delta\phi+3/4\eunorm{H}^2$.
\end{remark}

Closeness of the torsion implies two striking no--go theorems. The first one is this: if the dilaton is constant, then we get a classical solution, i.e. $H=0$.

\begin{theorem}\label{dilH}
For a spinor field $\Psi\in\Gamma(\Delta)$ satisfying~(\ref{gravdil}), it follows
$$
S^H=2\Delta\phi=-3\!\eunorm{H}^2. 
$$
In particular, $\phi=0$ if and only if $H=0$. Furthermore, if $M$ is compact, then $H=0$ so that any integrable generalised $G$--structure is classical.
\end{theorem}
\begin{proof}
Let $D^H$ denote the Dirac operator associated with $\nabla^H$, i.e.\ locally $D^H\Psi=\sum e_k\cdot\nabla^H_{e_k}\Psi$. By Theorem 3.3 in~\cite{friv02},
\begin{equation}\nonumber
D^H(H\cdot\Psi)=\big(d^*\!H\cdot-2\sigma^H\cdot-2\sum(e_k\llcorner H)\cdot\nabla^H_{e_k}\big)\Psi,
\end{equation}
where  $2\sigma^H=\sum(e_k\llcorner H)\wedge(e_k\llcorner H)$. Now
\begin{eqnarray*}
-2\sum (e_k\llcorner H)\cdot\nabla^H_{e_k}\cdot\Psi & = & \frac{1}{2}\sum(e_k\llcorner H)\cdot(e_k\llcorner H)\cdot\Psi\\
& = &(\frac{3}{2}\eunorm{H}^2+\sigma^H+\alpha)\cdot\Psi
\end{eqnarray*}
for some further $2$--form $\alpha$. On the other hand side, as in the computation of $\ric^H$,
\begin{eqnarray*}
D^H(H\cdot\Psi) & = & -2\sum e_k\cdot(\nabla^H_{e_k}d\phi)\cdot\Psi\\
& = & -2\sum (e_k\wedge\nabla^H_{e_k}d\phi-e_k\llcorner\nabla^H_{e_k}d\phi)\cdot\Psi\\
& = & (\alpha'-2\Delta^g\phi)\cdot\Psi,
\end{eqnarray*}
for some $2$--form $\alpha'$. Hence contracting with $q(\Psi,\cdot)$ yields
\begin{equation}\nonumber
q\big(\Psi,-2(\Delta^g\phi)\cdot\Psi\big)=q\big(\Psi,(\frac{3}{2}\!\eunorm{H}^2-\sigma^H)\cdot\Psi\big)
\end{equation}
\big($q(\alpha^p\cdot\Psi,\Psi)$ is purely imaginary for $p\equiv2(4)$ and real for $p\equiv0(4)$\big). As $q(\Psi,2\sigma^H\cdot\Psi)=S^H$ by Corollary 3.2 in~\cite{friv02}, the previous discussion implies $3\eunorm{H}^2=-S^H$. Hence, if $M$ is compact, $H=0$, as follows from integration of $\eunorm{H}^2$. 
\end{proof}

\begin{remark}
There are compact examples of generalised $G_2$\index{generalised $G_2$--structure}-- and $Spin(7)$\index{generalised $Spin(7)$--structure}--structures satisfying~(\ref{spinorfield}) for non--closed $H$~\cite{wi06}.
\end{remark}

\end{document}